\documentclass{scrartcl}
\usepackage{preamble}

\begin{document}

\maketitle

\begin{abstract}
     \small{\textsc{Abstract.} We compute extension sheaves of abelian schemes and of the additive group by the multiplicative group in the fppf topology. Our main results include a generalized and streamlined proof of the Barsotti--Weil formula, the vanishing of $\underline{\operatorname{Ext}}^2(A,\mathbb{G}_m)$ for an abelian scheme $A$ over a general base, and a description of $\underline{\operatorname{Ext}}^1(\mathbb{G}_a,\mathbb{G}_m)$ in characteristic zero.}
\end{abstract}

\section{Introduction}
Let $\mathcal{G}$ be a commutative group stack over a base scheme $S$. These objects, first studied by Deligne under the name \emph{champs de Picard strictement commutatifs} in \cite[Exp.~XVIII]{SGA43}, admit a Cartier dual defined by $\underline{\operatorname{Hom}}(\mathcal{G},\mathsf{B}\mathbb{G}_m)$. This construction generalizes the classical duals of abelian schemes, tori, Deligne's 1-motives, and numerous other group objects arising in algebraic geometry \cite{brochard2021duality}. Notably, this duality theory enabled Laumon to establish the geometric Langlands correspondence for tori \cite{Lau}.

Consider a morphism of abelian sheaves $\mathscr{A} \to \mathscr{B}$ in the fppf site $(\textsf{Sch}/S)_\text{fppf}$, viewed as an object of the derived category of abelian sheaves concentrated in degrees $-1$ and $0$. The sheaf $\mathscr{A}$ acts on $\mathscr{B}$ by translation, giving rise to the quotient stack $[\mathscr{B}/\mathscr{A}]$, which is a commutative group stack over $S$. By the Dold–Kan correspondence, every commutative group stack is equivalent to one arising in this way.

Under this correspondence, the Cartier dual of $[\mathscr{B}/\mathscr{A}]$ is identified with the commutative group stack associated to the complex
\[\tau_{\leq 0}\mathsf{R}\underline{\operatorname{Hom}}([\mathscr{A}\to \mathscr{B}],\mathbb{G}_m[1]).\]
As a result, computing Cartier duals of commutative group stacks reduces to computing extension sheaves of certain abelian sheaves by the multiplicative group. The purpose of this paper is to carry out such computations.

\subsection*{Extensions of abelian schemes}

The fact that the Cartier dual of an abelian scheme $A$ is its dual abelian scheme $A^\vee$ hinges on the isomorphism $\underline{\operatorname{Ext}}^1(A,\mathbb{G}_m)\simeq A^\vee$, known as the Barsotti--Weil formula. Despite its foundational role in arithmetic geometry, we were surprised to find a complete proof only in \cite[Exp.\ VII, \S1.3.8.2]{raynaud2006groupes}.\footnote{L.~Moret-Bailly later kindly informed us that another proof appears in the appendix of \cite{MoretBailly1981}. See also \cite[Footnote to Thm.~1.2.2]{jossen2009arithmetic} and \cite{BW} for further commentary.} Our first result generalizes this classical formula, and admits a simpler proof.

\begin{introthm}[\ref{generalized BW}]\label{thm A}
    Let $\mathscr{G}$ and $\mathscr{A}$ be abelian sheaves in $\normalfont(\textsf{Sch}/S)_\text{fppf}$. Assume that for $n = 1, 2, 3$, every morphism  of sheaves of sets $\mathscr{G}^n \to \mathscr{A}$ is constant, and that this remains true after any base change. Then, for any $S$-scheme $T$, the natural maps
    \[\underline{\operatorname{Ext}}^1(\mathscr{G},\mathscr{A})(T)\leftarrow\operatorname{Ext}^1_T(\mathscr{G},\mathscr{A})\to \normalfont\textrm{H}^1_m(\mathscr{G}_T,\mathscr{A}_T)\]
    are isomorphisms. 
\end{introthm}

We briefly explain the morphisms involved. The extension sheaf $\underline{\operatorname{Ext}}^1(\mathscr{G}, \mathscr{A})$ is the sheafification of the functor that assigns to each $S$-scheme $T$ the group $\operatorname{Ext}^1_T(\mathscr{G}, \mathscr{A})$ of extensions of $\mathscr{G}_T$ by $\mathscr{A}_T$. This gives rise to a natural map $\operatorname{Ext}^1_T(\mathscr{G},\mathscr{A})\to \underline{\operatorname{Ext}}^1(\mathscr{G},\mathscr{A})(T)$. Moreover, any such extension defines an $\mathscr{A}_T$-torsor over $\mathscr{G}_T$, yielding a map to $\normalfont\textrm{H}^1(\mathscr{G}_T,\mathscr{A}_T)$. Its image lies in the subgroup
\[\textrm{H}^1_m(\mathscr{G}_T,\mathscr{A}_T)\colonequals \ker\mleft(m^*-\operatorname{pr}_1^*-\operatorname{pr}_2^*\colon \textrm{H}^1(\mathscr{G}_T,\mathscr{A}_T)\to \textrm{H}^1(\mathscr{G}^2_T,\mathscr{A}_T)\mright),\]
where $m$, $\operatorname{pr}_1$, and $\operatorname{pr}_2$ denote the group law and the projections $\mathscr{G}_T \times_T \mathscr{G}_T \to \mathscr{G}_T$.

More generally, the theorem remains valid in any topos. The conclusion applies, in particular, when $\mathscr{G}$ is represented by an abelian scheme over $S$, and $\mathscr{A}$ is represented either by an affine commutative group scheme over $S$ or by a quasi-coherent $\mathcal{O}_S$-module.

\begin{introcor}[\ref{unipotent BW}, \ref{usual BW}]
	Let $A$ be an abelian scheme over $S$, and let $M$ be a quasi-coherent $\mathcal{O}_S$-module. Then there are isomorphisms
\[\underline{\operatorname{Ext}}^1(A,\mathbb{G}_m)\simeq A^\vee \qquad \text{and}\qquad \underline{\operatorname{Ext}}^1(A,M)\simeq \operatorname{Lie}(A^\vee)\otimes_{\mathcal{O}_S}M,\]
where $A^\vee$ denotes the dual abelian scheme of $A$, and $\operatorname{Lie}(A^\vee)$ its Lie algebra.
\end{introcor}

\subsection*{Higher extensions of abelian schemes}

In \cite{breen1969extensions}, Breen proved that for an abelian scheme $A$ over a \emph{regular} base $S$, the group $\operatorname{Ext}^2_S(A, \mathbb{G}_m)$ is torsion. This has frequently been interpreted as implying that the sheaf $\underline{\operatorname{Ext}}^2(A, \mathbb{G}_m)$ is torsion---see, for instance, \cite[Rem.~6]{breen1975theoreme}, \cite[Lem.~A.4.2]{chen2017geometric} and \cite[Cor.~11.5]{brochard2021duality}---and thus vanishes by a standard argument. However, this reasoning is flawed: even over a field, the fppf site\footnote{This happens even in the \emph{small} fppf site.} contains singular schemes. We provide a correct, albeit more involved, argument that establishes the vanishing of this sheaf over a general base.

\begin{introthm}[\ref{vanishing of extensions by additive group}, \ref{higher vanishing for abelian schemes}]\label{thm B}
	Let $A$ be an abelian scheme over a general base scheme $S$, and let $M$ be a quasi-coherent $\mathcal{O}_S$-module. Then the abelian sheaves
\[\underline{\operatorname{Ext}}^2(A,\mathbb{G}_m) \qquad \text{and}\qquad \underline{\operatorname{Ext}}^2(A,M)\]
both vanish.
\end{introthm}

In contrast to the previous corollary, for an $S$-scheme $T$, the sheafification map $\operatorname{Ext}^2_T(A, \mathbb{G}_m) \to \underline{\operatorname{Ext}}^2(A, \mathbb{G}_m)(T)$, and similarly for $M$, may fail to be an isomorphism. In Section~\ref{higher exts}, we compute these extension groups explicitly, and the resulting formulas imply the theorem.

It is worth noting that the sheaf $\underline{\operatorname{Ext}}^3(A, \mathbb{G}_m)$ does not always vanish. For example, using \cite{breenexample}, one can show that it is nonzero when $A$ is a supersingular elliptic curve over a separably closed field of characteristic two. Nevertheless, we expect the sheaves $\underline{\operatorname{Ext}}^i(A, \mathbb{G}_m)$ to be torsion for all $i \geq 2$. Under this assumption, \cite{breen1975theoreme} implies that they vanish in a range depending on the residual characteristics of the base scheme.

\subsection*{Extensions of the additive group}

The additive group $\mathbb{G}_a$ over $\mathbb{Q}$, being the simplest example of an algebraic group, naturally leads to the question of determining its Cartier dual in the sense of commutative group stacks. For any commutative group stack $\mathcal{G}$ with Cartier dual $\mathcal{G}^\vee$, there exists a Fourier transform
\[\textsf{D}_\text{qc}(\mathcal{G}) \to \textsf{D}_\text{qc}(\mathcal{G}^\vee),\]
which is often an equivalence. 

Since the sheaf $\underline{\operatorname{Hom}}(\mathbb{G}_a, \mathbb{G}_m)$ is represented by the formal completion $\widehat{\mathbb{G}}_a$ of $\mathbb{G}_a$ along the zero section, and since there is an equivalence of categories 
\[\textsf{D}_\text{qc}(\mathbb{G}_a) \simeq \textsf{D}_\text{qc}(\textsf{B}\widehat{\mathbb{G}}_a),\]
as shown in \cite[Ex.\ 2.2.12]{bhatt2022prismatic}, one might expect the Cartier dual of $\mathbb{G}_a$ to be the classifying stack $\textsf{B}\widehat{\mathbb{G}}_a$. This expectation holds if and only if the sheaf $\underline{\operatorname{Ext}}^1(\mathbb{G}_a, \mathbb{G}_m)$ vanishes.

An even more striking instance arises in the theory of $\mathcal{D}$-modules. Much of the recent progress in understanding irregular holonomic $\mathcal{D}$-modules relies on a Fourier transform
\[\textsf{D}_\text{qc}(\mathcal{D}_{\mathbb{A}^1})\to \textsf{D}_\text{qc}(\mathcal{D}_{\mathbb{A}^1}).\]
The derived category of quasi-coherent $\mathcal{D}$-modules on $\mathbb{A}^1$, denoted $\textsf{D}_\text{qc}(\mathcal{D}_{\mathbb{A}^1})$, is equivalent to $\textsf{D}_\text{qc}([\mathbb{G}_a /\widehat{\mathbb{G}}_a])$. This suggests that the commutative group stack $[\mathbb{G}_a / \widehat{\mathbb{G}}_a]$ should be self-dual. As before, this is equivalent to the vanishing of the sheaf $\underline{\operatorname{Ext}}^1(\mathbb{G}_a, \mathbb{G}_m)$.

Such a vanishing result appears in the published literature, notably in \cite[Lem.~1.3.6]{polishchuk2011kernel}, and implicitly in the proofs of \cite[Lem.~A.4.5]{barbieri2009sharp} and \cite[Lem.~10]{bertapelle2014generalized}. Yet, once again, the situation turns out to be more subtle: in \cite[Rem.~2.2.16]{rosengarten2023tate}, the second author describes a construction of Gabber yielding a nonzero section of this sheaf. Our final theorem provides a complete computation of this object.

\begin{introthm}[\ref{Ext additive}]\label{thm C}
Let $T$ be a quasi-compact and quasi-separated $\mathbb{Q}$-scheme. Then the natural maps
\[\normalfont\underline{\operatorname{Ext}}^1(\mathbb{G}_a,\mathbb{G}_m)(T)\to \underline{\operatorname{Ext}}^1(\mathbb{G}_a,\mathbb{G}_m)(T_\text{red})\leftarrow \operatorname{Ext}^1_{T_\text{red}}(\mathbb{G}_a,\mathbb{G}_m)\to \textrm{H}^1_m(\mathbb{G}_{a,T_\text{red}},\mathbb{G}_{m})\]
are all isomorphisms. These groups vanish if $\normalfont T_\text{red}$ is seminormal. Conversely, if $\normalfont T_\text{red}$ is affine and not seminormal, they are nonzero.
\end{introthm}

A similar result was announced by Gabber in a recent conference talk \cite{gabberstalk}. Our proof was developed independently, and we were not aware of Gabber's argument at the time. Most of the content of this theorem also appears in the first author's thesis \cite{ribeiro}. We also note that there exists an example due to Weibel of a reduced quasi-compact and quasi-separated $\mathbb{Q}$-scheme $T$ that is not affine and not seminormal, for which $\textrm{H}^1_m(\mathbb{G}_{a,T},\mathbb{G}_{m})$ vanishes; see Remark~\ref{weibels example}.

\subsection*{Outline of the main arguments}

Let $\mathscr{G}$ and $\mathscr{A}$ be abelian sheaves on the big fppf site $(\textsf{Sch}/S)_\text{fppf}$, and let $T$ be an $S$-scheme. One of the main tools used in this paper is a pair of spectral sequences:
\[E_2^{i,j}=\normalfont \textrm{H}^i(T,\underline{\operatorname{Ext}}^j(\mathscr{G},\mathscr{A})) \quad \text{and}\quad F_1^{i,j} = \prod_{r=1}^{n_i}\normalfont\textrm{H}^j(\mathscr{G}^{s_{i,r}}_T,\mathscr{A}_T),\]
both converging to $\operatorname{Ext}^{i+j}_T(\mathscr{G},\mathscr{A})$. The first is the classical local-to-global spectral sequence. The second arises from the Breen--Deligne resolution, whose existence was independently established by Deligne---in a letter to Breen made public only recently in \cite[App.~B]{ribeiro}---and by Clausen--Scholze in \cite[Thm.~4.10]{scholze2019condensed}.

The morphisms in the statements of Theorems~\ref{thm A} and \ref{thm C} appear in the exact sequences of low-degree terms associated with these spectral sequences, which take the form:
\[\begin{tikzcd}[ampersand replacement=\&,sep=small]
	0 \& {\textrm{H}^1(T,\underline{\operatorname{Hom}}(\mathscr{G},\mathscr{A}))} \& {\operatorname{Ext}^1_T(\mathscr{G},\mathscr{A})} \& {\underline{\operatorname{Ext}}^1(\mathscr{G},\mathscr{A})(T)} \& {\textrm{H}^2(T,\underline{\operatorname{Hom}}(\mathscr{G},\mathscr{A}))} \\
	0 \& {\textrm{H}^2_s(\mathscr{G}_T,\mathscr{A}_T)} \& {\operatorname{Ext}^1_T(\mathscr{G},\mathscr{A})} \& {\textrm{H}^1_m(\mathscr{G}_T,\mathscr{A}_T)} \& {\textrm{H}^3_s(\mathscr{G}_T,\mathscr{A}_T).}
	\arrow[from=1-1, to=1-2]
	\arrow[from=1-2, to=1-3]
	\arrow[from=1-3, to=1-4,color=ocre]
	\arrow[from=1-4, to=1-5]
	\arrow[from=2-1, to=2-2]
	\arrow[from=2-2, to=2-3]
	\arrow[from=2-3, to=2-4,color=ocre]
	\arrow[from=2-4, to=2-5]
\end{tikzcd}\]
Here, the groups $\textrm{H}^2_s(\mathscr{G}_T,\mathscr{A}_T)$ and $\textrm{H}^3_s(\mathscr{G}_T,\mathscr{A}_T)$ are analogues of group cohomology, defined as the cohomology of a complex whose terms are of the form $\operatorname{Mor}_T(\mathscr{G}_T^n,\mathscr{A}_T)$ for suitable integers $n$. (Definition~\ref{def invariants}.)

In the setting of Theorem~\ref{thm A}, we assume that all necessary morphisms of sheaves of sets $\mathscr{G}_T^n \to \mathscr{A}_T$ are constant. Under this hypothesis, the sheaf $\underline{\operatorname{Hom}}(\mathscr{G}, \mathscr{A})$ vanishes, and the complex computing $\mathrm{H}^2_s(\mathscr{G}_T, \mathscr{A}_T)$ and $\mathrm{H}^3_s(\mathscr{G}_T, \mathscr{A}_T)$ simplifies considerably. It follows that the relevant cohomology groups above vanish, thereby establishing Theorem~\ref{thm A}. The proof of Theorem~\ref{thm C} proceeds in a similar spirit, though the necessary computations are substantially more involved.

The proof of Theorem~\ref{thm B} in the case of the multiplicative group proceeds by dévissage, starting with the case in which the base scheme is the spectrum of a field. This method has become fairly standard for problems of this type: one first extends the result from fields to artinian local rings, then to their completions, and finally deduces the statement over the base scheme $S$.

We now sketch the argument, omitting many technical details. Let $\mathscr{E} \in \operatorname{Ext}^2_T(A, \mathbb{G}_m)$ be an extension class. By a limit argument, we may assume that $T$ is excellent. Since, for every non-zero integer $n$, the multiplication-by-$n$ map on $\underline{\operatorname{Ext}}^2(A,\mathbb{G}_m)$ is injective (Lemma~\ref{Ext is torsion-free}), it suffices to show that for each point $t\in T$, there exists a flat morphism of finite presentation $g\colon T^\prime\to T$, whose image contains $t$, such that $g^*\mathscr{E}$ is torsion. 

To construct such a morphism, let $(B, \mathfrak{m})$ denote the local ring of $T$ at $t$. By Breen's result over regular schemes, the restriction of $\mathscr{E}$ to the residue field $B/\mathfrak{m}$ is torsion. We then prove that, for any noetherian ring $R$, the natural map
\[\operatorname{Ext}^2_{R}(A, \mathbb{G}_m) \to \operatorname{Ext}^2_{{R}_{\mathrm{red}}}(A, \mathbb{G}_m)\]
is injective. This follows from the computation of higher extension groups of abelian schemes by quasi-coherent sheaves, which are related---via one of the spectral sequences discussed above---to quasi-coherent cohomology of abelian schemes.

The injectivity result above implies that $\mathscr{E}$ restricts to a torsion class over $B/\mathfrak{m}^n$ for all $n > 0$, since $(B/\mathfrak{m}^n)_{\mathrm{red}} \simeq B/\mathfrak{m}$. Let $\widehat{B}$ denote the $\mathfrak{m}$-adic completion of $B$. We then show that the natural map
\[\operatorname{Ext}^2_{\widehat{B}}(A, \mathbb{G}_m) \to \lim_n \operatorname{Ext}^2_{B/\mathfrak{m}^n}(A, \mathbb{G}_m)\]
is injective, and hence that $\mathscr{E}$ becomes torsion after pullback to $\widehat{B}$.

Unlike the previous injectivity statement, the cohomological analogue of this map---obtained by replacing $\operatorname{Ext}^2$ with $\textrm{H}^2$---need not be injective \cite[Ex.~9.3]{kresch2023formal}. In Section~\ref{section on completions}, we adapt the methods of \cite{kresch2023formal}, which rely on an algebraization theorem for algebraic stacks due to Bhatt and Halpern-Leistner \cite[Thm.~7.4]{bhl}. The necessary modifications consist mainly in replacing the small étale site with the big fppf site in the formalism of continuous cohomology.

Finally, by spreading out the morphism $\operatorname{Spec}\widehat{B} \to \operatorname{Spec}B$, we obtain an affine open neighborhood $V$ of $t$ and a smooth morphism $Y \to V$ whose image contains $t$, such that the restriction of $\mathscr{E}$ to $Y$ is torsion. This concludes the proof.

\vspace{1em}

\textsc{Acknowledgements.} Most of this work was carried out while the first author was a doctoral student of Javier Fresán, whose guidance and support are gratefully acknowledged. We thank Ofer Gabber, Andrew Kresch, and Siddharth Mathur for helpful discussions related to Lemma~\ref{extensions and completions}. The first author was supported by Swiss National Science Foundation grant 219220, and the second by Israel Science Foundation grant 2083/24.

\section{The fundamental spectral sequences}\label{Sect 2}
Let $\textsf{X}$ be a topos, and denote by $\textsf{Ab}(\textsf{X})$ the category of abelian group objects in $\textsf{X}$. Given an object $T$ of $\textsf{X}$ and an abelian group $\mathscr{A}$ in $\textsf{X}$, the cohomology group $\textrm{H}^i(T,\mathscr{A})$ is defined as the $i$-th right derived functor of $\operatorname{Mor}_\textsf{X}(T,-)\colon \textsf{Ab}(\textsf{X})\to \textsf{Ab}$, evaluated at $\mathscr{A}$. As usual, we write $\Gamma(T,\mathscr{A})$ for the group $\textrm{H}^0(T,\mathscr{A})=\operatorname{Mor}_\textsf{X}(T,\mathscr{A})$.

\begin{remark}\label{cohomology is reasonable}
    Let $S$ be a scheme, and let $\textsf{X}$ be the category of sheaves on the site $(\textsf{Sch}/S)_\tau$, where $\tau$ is a given Grothendieck topology. By the Yoneda lemma, $\operatorname{Mor}_\textsf{X}(S,-)$ agrees with the usual global sections functor, thereby justifying the above definition. More generally, for a morphism $f\colon T\to S$, the scheme $T$ defines an object of $\textsf{X}$, and the derived functor $\textrm{H}^i(T,-)\colon \textsf{Ab}((\textsf{Sch}/S)_\tau)\to \textsf{Ab}$ coincides with the composition
    \[\textsf{Ab}((\textsf{Sch}/S)_\tau)\xrightarrow{\ f^*\ }\textsf{Ab}((\textsf{Sch}/T)_\tau)\xrightarrow{\ \textrm{H}^i(T,-)\ }\textsf{Ab},\]
    as shown in \cite[Tag \href{https://stacks.math.columbia.edu/tag/03F3}{03F3}]{Stacks}.
\end{remark}

We denote by $\mathscr{A}_T$ the product $\mathscr{A}\times T$ considered as an abelian group object in the localized topos $\textsf{X}/T$ \cite[Tag \href{https://stacks.math.columbia.edu/tag/04GY}{04GY}]{Stacks}. For another abelian group $\mathscr{G}$ in $\textsf{X}$, we denote by $\operatorname{Hom}_T(\mathscr{G},\mathscr{A})$ the group of morphisms $\mathscr{G}_T\to \mathscr{A}_T$ in $\textsf{Ab}(\textsf{X}/T)$. The assignment $T\mapsto \operatorname{Hom}_T(\mathscr{G},\mathscr{A})$ defines a sheaf\footnote{According to \cite[Prop.\ IV.1.4]{SGA41}, this is the same as the functor $T\mapsto \operatorname{Hom}_T(\mathscr{G},\mathscr{A})$ sending colimits in $\textsf{X}$ to limits in $\textsf{Ab}$.} on $\textsf{X}$ with respect to the canonical topology, denoted by $\underline{\operatorname{Hom}}(\mathscr{G},\mathscr{A})$. Under the equivalence between $\textsf{X}$ and its category of sheaves, $\underline{\operatorname{Hom}}(\mathscr{G},\mathscr{A})$ defines an object of $\textsf{Ab}(\textsf{X})$.

\begin{definition}\label{def ext}
Let $T$ be an object of a topos $\normalfont\textsf{X}$, and let $\mathscr{G}$ be an abelian group in $\textsf{X}$. The functors
\[
\operatorname{Ext}_T^i(\mathscr{G}, -)  \colon \textsf{Ab}(\textsf{X}) \to \textsf{Ab} \quad \text{and} \quad \underline{\operatorname{Ext}}^i(\mathscr{G}, -) \colon \textsf{Ab}(\textsf{X}) \to \textsf{Ab}(\textsf{X})
\]
are defined as the $i$-th right derived functors of $\operatorname{Hom}_T(\mathscr{G}, -)$ and $\underline{\operatorname{Hom}}(\mathscr{G}, -)$, respectively. When $T$ is the final object of $\textsf{X}$, we write $\operatorname{Ext}_T^i(\mathscr{G}, -)$ simply as $\operatorname{Ext}^i(\mathscr{G}, -)$.
\end{definition}

We will often study the sheaf $\underline{\operatorname{Ext}}^i(\mathscr{G}, \mathscr{A})$ using the fact that it is the sheafification of the presheaf $T \to \operatorname{Ext}_T^i(\mathscr{G}, \mathscr{A})$ \cite[Prop.\ V.6.1]{SGA41}. In particular, there exists a \emph{sheafification map}
\[
\operatorname{Ext}_T^i(\mathscr{G}, \mathscr{A}) \to \underline{\operatorname{Ext}}^i(\mathscr{G}, \mathscr{A})(T),
\]
which is functorial in $\mathscr{G}$, $\mathscr{A}$, and $T$. Combined with \cite[Tag \href{https://stacks.math.columbia.edu/tag/00WK}{00WK}]{Stacks}, this yields a rather concrete description of $\underline{\operatorname{Ext}}^i(\mathscr{G}, \mathscr{A})(T)$. The following result, also contained in \cite[Prop.\ V.6.1]{SGA41}, gives another relation between these objects.

\begin{proposition}[Local-to-global spectral sequence]\label{sheafification map}
	Let $T$ be an object and $\mathscr{G},\mathscr{A}$ be abelian groups in a topos $\normalfont\textsf{X}$. There exists a spectral sequence
\[E_2^{i,j}=\normalfont \textrm{H}^i(T,\underline{\operatorname{Ext}}^j(\mathscr{G},\mathscr{A}))\implies \operatorname{Ext}_T^{i+j}(\mathscr{G},\mathscr{A}),\]
that is functorial in $\mathscr{G}$, $\mathscr{A}$, and $T$. In particular, there is an exact sequence of low-degree terms
\[\begin{tikzcd}[row sep=scriptsize]
	0 & {\normalfont\textrm{H}^1(T,\underline{\operatorname{Hom}}(\mathscr{G},\mathscr{A}))} & {\operatorname{Ext}^1_T(\mathscr{G},\mathscr{A})\xrightarrow{\hspace{2.68em}} \underline{\operatorname{Ext}}^1(\mathscr{G},\mathscr{A})(T)}\ar[draw=none]{dl}[name=X, anchor=center]{}\ar[rounded corners,
            to path={ -- ([xshift=2ex]\tikztostart.east)
                      |- (X.center) \tikztonodes
                      -| ([xshift=-2ex]\tikztotarget.west)
                      -- (\tikztotarget)}]{dl}[at end]{} \\
	& {\normalfont\textrm{H}^2(T,\underline{\operatorname{Hom}}(\mathscr{G},\mathscr{A}))} & {\ker\mleft(\operatorname{Ext}^2_T(\mathscr{G},\mathscr{A})\to \underline{\operatorname{Ext}}^2(\mathscr{G},\mathscr{A})(T)\mright)}\ar[draw=none]{dl}[name=X, anchor=center]{}\ar[rounded corners,
            to path={ -- ([xshift=2ex]\tikztostart.east)
                      |- (X.center) \tikztonodes
                      -| ([xshift=-2ex]\tikztotarget.west)
                      -- (\tikztotarget)}]{dl}[at end]{} \\
	& {\normalfont\textrm{H}^1(T,\underline{\operatorname{Ext}}^1(\mathscr{G},\mathscr{A}))} & {\normalfont\textrm{H}^3(T,\underline{\operatorname{Hom}}(\mathscr{G},\mathscr{A})),}
	\arrow[from=1-1, to=1-2]
	\arrow[from=1-2, to=1-3]
	\arrow[from=2-2, to=2-3]
	\arrow[from=3-2, to=3-3]
\end{tikzcd}\]
that is functorial in $\mathscr{G}$, $\mathscr{A}$, and $T$.
\end{proposition}

Now that we have established a connection between extension \emph{sheaves} and extension \emph{groups}, we will study a method for computing the latter. The following proposition, suggested by Grothendieck in \cite[Exp.\ VII, Rem.\ 3.5.4]{raynaud2006groupes} and partially developed by Breen in \cite{breen1969extensions}, has been independently proven by Deligne and by Clausen--Scholze \cite[Thm.\ 4.10]{scholze2019condensed}.

\begin{proposition}[Breen--Deligne resolution]\label{BD resolution}
Let $\mathscr{G}$ be an abelian group in a topos $\normalfont\textsf{X}$. There exists a functorial resolution of the form
\[\cdots \to \bigoplus_{j=1}^{n_i}\mathbb{Z}[\mathscr{G}^{s_{i,j}}]\to \cdots \to \mathbb{Z}[\mathscr{G}^3]\oplus \mathbb{Z}[\mathscr{G}^2]\to \mathbb{Z}[\mathscr{G}^2]\to \mathbb{Z}[\mathscr{G}]\to \mathscr{G},\]
where the $n_i$ and $s_{i,j}$ are all positive integers.
\end{proposition}

More precisely, Clausen--Scholze's argument shows that any given partial resolution of $\mathscr{G}$ can be extended to a resolution of the form described above. In particular, the initial terms of the resolution may be taken to be
\[\mathbb{Z}[\mathscr{G}^4]\oplus \mathbb{Z}[\mathscr{G}^3]\oplus \mathbb{Z}[\mathscr{G}^3]\oplus \mathbb{Z}[\mathscr{G}^2]\oplus \mathbb{Z}[\mathscr{G}]\xrightarrow{d_3} \mathbb{Z}[\mathscr{G}^3]\oplus \mathbb{Z}[\mathscr{G}^2]\xrightarrow{d_2} \mathbb{Z}[\mathscr{G}^2]\xrightarrow{d_1} \mathbb{Z}[\mathscr{G}]\xrightarrow{d_0} \mathscr{G},\]
with differentials explicitly given by:
\begin{align*}
    d_3([x,y,z,t]) &= ([x+y,z,t] - [x,y+z,t] + [x,y,z+t]- [x,y,z] - [y,z,t],0)\\
    d_3([x,y,z]) &= (-[x,y,z]+[x,z,y]-[z,x,y],[x+y,z]-[x,z]-[y,z])\\
    d_3([x,y,z]) &= ([x,y,z]-[y,x,z]+[y,z,x],[x,y+z]-[x,y]-[x,z])\\
    d_3([x,y]) &= (0,[x,y]+[y,x])\\
    d_3([x]) &= (0,[x,x]) \\
    d_2([x,y,z]) &= [x+y,z]-[x,y+z]+[x,y]-[y,z] \\
	d_2([x,y]) &= [x,y]-[y,x] \\
    d_1([x,y]) &= [x+y]-[x]-[y]\\
	d_0([x]) &= x.
\end{align*}
Here, the top $d_3([x,y,z])$ acts on the first factor of $\mathbb{Z}[\mathscr{G}^3]$, while the bottom $d_3([x,y,z])$ acts on the second factor. Throughout this paper, any Breen--Deligne resolution is assumed to begin with these terms. This explicit presentation enables the definition of two invariants.

\begin{definition}\label{def invariants}
	Let $\mathscr{G}$ and $\mathscr{A}$ be abelian groups in a topos $\textsf{X}$. Applying the functor $\operatorname{Hom}(-,\mathscr{A})$ to a Breen--Deligne resolution of $\mathscr{G}$, we obtain the complex
\[\begin{tikzcd}
	{\Gamma(\mathscr{G},\mathscr{A})\to  \Gamma(\mathscr{G}^2,\mathscr{A}) \to \Gamma(\mathscr{G}^3,\mathscr{A})\oplus \Gamma(\mathscr{G}^2,\mathscr{A})} \ar[draw=none]{d}[name=X, anchor=center]{}\ar[rounded corners,
            to path={ -- ([xshift=2ex]\tikztostart.east)
                      |- (X.center) \tikztonodes
                      -| ([xshift=-2ex]\tikztotarget.west)
                      -- (\tikztotarget)}]{d}[at end]{} \\
	{\Gamma(\mathscr{G}^4,\mathscr{A})\oplus \Gamma(\mathscr{G}^3,\mathscr{A})\oplus \Gamma(\mathscr{G}^3,\mathscr{A})\oplus \Gamma(\mathscr{G}^2,\mathscr{A})\oplus \Gamma(\mathscr{G},\mathscr{A}),}
\end{tikzcd}\]
concentrated in degrees 1 to 4. The first and second cohomology groups of this complex are denoted by $\textrm{H}^2_s(\mathscr{G},\mathscr{A})$ and $\textrm{H}^3_s(\mathscr{G},\mathscr{A})$, respectively.
\end{definition}

\begin{remark}
	The invariant $\textrm{H}^2_s(\mathscr{G},\mathscr{A})$ is usually known as the symmetric subgroup of the second Hochschild cohomology\footnote{Defined, akin to the bar resolution in group cohomology, as the cohomology of the simpler complex $\Gamma(\mathscr{G},\mathscr{A})\to  \Gamma(\mathscr{G}^2,\mathscr{A}) \to \Gamma(\mathscr{G}^3,\mathscr{A})$.} group $\textrm{H}^2_0(\mathscr{G},\mathscr{A})$ \cite[Chap.\ 15]{M}. When $\textsf{X}$ is the topos of sets, it reduces to the subgroup of the group cohomology $\textrm{H}^2(\mathscr{G},\mathscr{A})$ constituted of the symmetric cocycles. The notation $\textrm{H}^3_s(\mathscr{G},\mathscr{A})$ indicates that this group is, in some sense, a variant of the third Hochschild cohomology which is more adapted to commutative groups.
\end{remark}

As usual, the group $\textrm{H}^1(T,\mathscr{A})$ classifies $\mathscr{A}$-torsors over $T$, where the group operation corresponds to the contracted product \cite[\S\S III.2.4, III.3.5]{giraud2020cohomologie}. For a morphism $f\colon T \to S$ in $\textsf{X}$, the induced map $f^*\colon \textrm{H}^1(S,\mathscr{A})\to \textrm{H}^1(T,\mathscr{A})$ sends an $\mathscr{A}$-torsor $P\to S$ to the pullback $f^*P\to T$ \cite[\S V.1.5]{giraud2020cohomologie}. When $T=\mathscr{G}$ is also an abelian group, we define a subgroup $\textrm{H}_m^1(\mathscr{G},\mathscr{A})$ of $\textrm{H}^1(\mathscr{G},\mathscr{A})$ constituted of the $\mathscr{A}$-torsors over $\mathscr{G}$ compatible with the group structure on the latter.

\begin{definition}\label{def Hm}
	Let $\mathscr{G}$ and $\mathscr{A}$ be abelian groups in a topos $\normalfont\textsf{X}$. Denote by $m\colon\mathscr{G}\times \mathscr{G}\to \mathscr{G}$ the group operation of $\mathscr{G}$, and by $\operatorname{pr}_1,\operatorname{pr}_2\colon \mathscr{G}\times \mathscr{G}\to \mathscr{G}$ the natural projections. We define $\textrm{H}_m^1(\mathscr{G},\mathscr{A})$ as the kernel of the morphism $m^*-\operatorname{pr}_1^*-\operatorname{pr}_2^*$.
\end{definition}

Put simply, $\textrm{H}^1_m(\mathscr{G},\mathscr{A})$ is the group of isomorphism classes of $\mathscr{A}$-torsors $P$ over $\mathscr{G}$ satisfying $m^*P\simeq\operatorname{pr}_1^*P\wedge\operatorname{pr}_2^*P$. These $\mathscr{A}$-torsors are often referred to in the literature as being \emph{multiplicative} or \emph{primitive}. With this terminology established, we may now explain the computation of the extension groups.

\begin{proposition}\label{short exact sequence computing extensions}
	Let $\mathscr{G}$ and $\mathscr{A}$ be abelian groups in a topos $\normalfont\textsf{X}$. There exists a spectral sequence
    \[E_1^{i,j} = \prod_{r=1}^{n_i}\normalfont\textrm{H}^j(\mathscr{G}^{s_{i,r}},\mathscr{A})\implies \operatorname{Ext}^{i+j}(\mathscr{G},\mathscr{A}),\]
    where $n_i$ and $s_{i,r}$ are the positive integers appearing in a Breen--Deligne resolution, that is functorial in $\mathscr{G}$ and $\mathscr{A}$. In particular, there is an exact sequence of low-degree terms
\[0\to \normalfont \textrm{H}^2_s(\mathscr{G},\mathscr{A})\to \operatorname{Ext}^1(\mathscr{G},\mathscr{A})\to \textrm{H}^1_m(\mathscr{G},\mathscr{A})\to \textrm{H}^3_s(\mathscr{G},\mathscr{A})\to\operatorname{Ext}^2(\mathscr{G},\mathscr{A}),\]
that is functorial in $\mathscr{G}$ and $\mathscr{A}$.
\end{proposition}

Before diving into the proof, let us explain the morphism $\operatorname{Ext}^1(\mathscr{G},\mathscr{A})\to \textrm{H}^1_m(\mathscr{G},\mathscr{A})$. Since $\operatorname{Ext}^1(-,\mathscr{A})$ is an additive functor, an extension $\mathscr{E}$ of $\mathscr{G}$ by $\mathscr{A}$ always satisfies
\[m^*\mathscr{E}=(\operatorname{pr}_1+\operatorname{pr}_2)^*\mathscr{E}\simeq \operatorname{pr}_1^*\mathscr{E}+\operatorname{pr}_2^*\mathscr{E},\]
where the sum on the right is the Baer sum of extensions. Such an extension defines an $\mathscr{A}$-torsor $P$ over $\mathscr{G}$, which satisfies $m^*P\simeq\operatorname{pr}_1^*P\wedge\operatorname{pr}_2^*P$.

\begin{proof}[Proof of Proposition~\ref{short exact sequence computing extensions}]
	The universal property of free objects gives that $\textrm{H}^i(\mathscr{G}^n,\mathscr{A})$ is isomorphic to $\operatorname{Ext}^i(\mathbb{Z}[\mathscr{G}^n],\mathscr{A})$ for all $n$ and $i$. The desired spectral sequence then arises as an instance of the hypercohomology spectral sequence \cite[Tag \href{https://stacks.math.columbia.edu/tag/07AA}{07AA}]{Stacks}.
\end{proof}

In the topos of sets, every torsor is trivial, and the proposition above recovers the classical isomorphism $\textrm{H}^2_s(\mathscr{G}, \mathscr{A}) \simeq \operatorname{Ext}^1(\mathscr{G}, \mathscr{A})$ for any abelian groups $\mathscr{G}$ and $\mathscr{A}$. In contrast, in the setting of this paper---where $\textsf{X}$ is a large fppf topos---the groups $\textrm{H}^2_s(\mathscr{G}, \mathscr{A})$ and $\textrm{H}^3_s(\mathscr{G}, \mathscr{A})$ are often trivial, while $\textrm{H}^1_m(\mathscr{G}, \mathscr{A})$ is the object of real interest.

\begin{remark}\label{normalized BD}
For abelian groups $\mathscr{G}$ and $\mathscr{A}$ in a topos $\textsf{X}$, we define the \emph{normalized} cohomology group $\textrm{H}^i_N(\mathscr{G},\mathscr{A})$ as the kernel of the morphism
\[e^*\colon \textrm{H}^i(\mathscr{G},\mathscr{A})\to \textrm{H}^i(0,\mathscr{A}),\]
where $0$ is the final object of $\textsf{X}$ and $e\colon 0\to \mathscr{G}$ is the zero-section of $\mathscr{G}$. Just as group cohomology can be computed by the normalized bar resolution, extension groups can be computed by a normalized variant of the Breen--Deligne spectral sequence:
\[E_1^{i,j}\colon\prod_{r=1}^{n_i}\textrm{H}^j_N(\mathscr{G}^{s_{i,r}},\mathscr{A})\implies \operatorname{Ext}^{i+j}(\mathscr{G},\mathscr{A}).\]

To prove this, denote by $B(\mathscr{G})$ a Breen--Deligne resolution of $\mathscr{G}$ (including the last term $\mathscr{G}$). From the long exact sequence in cohomology associated to the short exact sequence of complexes
\[0\to B(0)\to B(\mathscr{G})\to B(\mathscr{G})/B(0)\to 0,\]
we see that the complex
\[B(\mathscr{G})/B(0)=\left[ \cdots \to \bigoplus_{j=1}^{n_i}\frac{\mathbb{Z}[\mathscr{G}^{s_{i,j}}]}{\mathbb{Z}}\to \cdots \to \frac{\mathbb{Z}[\mathscr{G}^3]}{\mathbb{Z}}\oplus \frac{\mathbb{Z}[\mathscr{G}^2]}{\mathbb{Z}}\to \frac{\mathbb{Z}[\mathscr{G}^2]}{\mathbb{Z}}\to \frac{\mathbb{Z}[\mathscr{G}]}{\mathbb{Z}}\to \mathscr{G}\right]\]
is still a resolution of $\mathscr{G}$. Consequently, the same argument as in the proof of Proposition~\ref{short exact sequence computing extensions} gives a spectral sequence computing $\operatorname{Ext}^{i+j}(\mathscr{G},\mathscr{A})$ whose terms are products of copies of $\operatorname{Ext}^j(\mathbb{Z}[\mathscr{G}^{s_{i,r}}]/\mathbb{Z},\mathscr{A})$. These objects fit into an exact sequence
\[\begin{tikzcd}[row sep=scriptsize]
	& {\underbrace{\operatorname{Ext}^{j-1}(\mathbb{Z}[\mathscr{G}^{s_{i,r}}],\mathscr{A})}_{\textrm{H}^{j-1}(\mathscr{G}^{s_{i,r}},\mathscr{A})}} & {\underbrace{\operatorname{Ext}^{j-1}(\mathbb{Z},\mathscr{A})}_{\textrm{H}^{j-1}(0,\mathscr{A})}}  \ar[draw=none]{dll}[name=X, anchor=center]{}\ar[rounded corners,
            to path={ -- ([xshift=2ex]\tikztostart.east)
                      |- ([yshift=-1em] X.center) \tikztonodes
                      -| ([xshift=-2ex]\tikztotarget.west)
                      -- (\tikztotarget)}]{dll}[at end]{} \\
	{\operatorname{Ext}^j(\mathbb{Z}[\mathscr{G}^{s_{i,r}}]/\mathbb{Z},\mathscr{A})} & {\underbrace{\operatorname{Ext}^j(\mathbb{Z}[\mathscr{G}^{s_{i,r}}],\mathscr{A})}_{\textrm{H}^j(\mathscr{G}^{s_{i,r}},\mathscr{A})}} & {\underbrace{\operatorname{Ext}^j(\mathbb{Z},\mathscr{A})}_{\textrm{H}^j(0,\mathscr{A})}.}
	\arrow[from=1-2, to=1-3]
	\arrow[from=2-1, to=2-2]
	\arrow[from=2-2, to=2-3]
\end{tikzcd}\]
Since the map $e^*\colon \textrm{H}^{j-1}(\mathscr{G}^{s_{i,r}},\mathscr{A})\to \textrm{H}^{j-1}(0,\mathscr{A})$ has a section, it is surjective and then the exact sequence above shows that $\operatorname{Ext}^j(\mathbb{Z}[\mathscr{G}^{s_{i,r}}]/\mathbb{Z},\mathscr{A})$ is isomorphic to $\textrm{H}^j_N(\mathscr{G}^{s_{i,r}},\mathscr{A})$.
\end{remark}

Henceforth, we fix a base scheme $S$ and consider $\textsf{X}$ to be the category of sheaves on the large fppf site $(\textsf{Sch}/S)_\text{fppf}$ unless otherwise specified. The following corollary shows that we can often consider coarser topologies as well.

\begin{corollary}\label{different topologies}
Let $G$ be a commutative group scheme over $S$. For a smooth commutative group scheme $H$ over $S$ and an $S$-scheme $T$, the extension group $\operatorname{Ext}^i_T(G,H)$ can be computed in the étale topology as well. Moreover, for a quasi-coherent $\mathcal{O}_S$-module $M$, the extension group $\operatorname{Ext}^i_T(G,M)$ can be computed in the étale and Zariski topologies as well.
\end{corollary}

\begin{proof}
    Using the Breen--Deligne spectral sequence, those extension groups can be computed in terms of some cohomology groups, which are independent of the topology.
\end{proof}

Since our next discussion may involve potentially confusing notation, we will temporarily include the site in the notation for cohomology groups. Given a Grothendieck site $\textsf{C}$ and an object $T$ of the topos $\textsf{X} = \textsf{Sh}(\textsf{C})$, we denote the $i$-th right derived functor of $\operatorname{Mor}_\textsf{X}(T, -)\colon \textsf{Ab}(\textsf{X}) \to \textsf{Ab}$ by $\textrm{H}^i(\textsf{C}; T, -)$.

Let $f\colon T\to S$ be a morphism of schemes. Although $T$ may not itself be an object of the small étale site $S_\text{ét}$, it still defines a sheaf therein. Consequently, the construction above yields a functor $\textrm{H}^i(S_\text{ét};T,-)$. This functor does \emph{not} necessarily agree with the usual definition of $\textrm{H}^i(T,f^*-)$, which is $\textrm{H}^i(T_\text{ét};T,f^*-)$. Indeed, even the diagram
\[\begin{tikzcd}
	{\textsf{Ab}(S_\text{ét})} &&& {\textsf{Ab}} \\
	{\textsf{Ab}(T_\text{ét})} &&& {\textsf{Ab}}
	\arrow["{\operatorname{Mor}_{\textsf{Sh}(S_\text{ét})}(T,-)}", from=1-1, to=1-4]
	\arrow["{f^*}"', from=1-1, to=2-1]
	\arrow[equals, from=1-4, to=2-4]
	\arrow["{\Gamma(T,-)}", from=2-1, to=2-4]
\end{tikzcd}\]
might fail to commute unless $f$ is étale. This limits the utility of the small étale site in the context of this paper.

That said, since both abelian schemes and the additive group are smooth over the base, the \emph{lisse-étale} site proves to be more appropriate for our purposes.

\begin{corollary}\label{small fppf site}
Let $G$ be a commutative group scheme that is smooth over $S$. For an abelian sheaf $\mathscr{A}$ on $(\normalfont\textsf{Sch}/S)_\text{ét}$ and an $S$-scheme $T$, the extension group $\operatorname{Ext}^i_T(G,\mathscr{A})$ computed on this site agrees with the corresponding group computed in the lisse-étale site.
\end{corollary}

Given that the lisse-étale site may be unfamiliar to some readers---and differs in essential ways from the usual étale sites---we encourage the reader to consult the following remark.\footnote{See also \cite[\S3]{olsson} or \cite[Tags \href{https://stacks.math.columbia.edu/tag/0786}{0786} and \href{https://stacks.math.columbia.edu/tag/0GR1}{0GR1}]{Stacks} for further details.} It offers the context needed to understand the Lemma~\ref{small fppf cohomology}, after which the proof of Corollary~\ref{small fppf site} follows an argument analogous to that of Corollary~\ref{different topologies}.

\begin{remark}[The lisse-étale site]\label{lisse-etale remark}
    The \emph{lisse-étale site} of a scheme $S$, denoted $({\textsf{Sm}}/S)_{\text{ét}}$, is the full subcategory of $\textsf{Sch}/S$ consisting of smooth $S$-schemes. A covering $\{U_i \to U\}_{i \in I}$ in this site is a family of étale morphisms such that the induced map
    \[
    \coprod_{i \in I} U_i \to U
    \]
    is surjective.

    The inclusion functor $\textsf{Sm}/S \hookrightarrow \textsf{Sch}/S$ is both continuous and cocontinuous. It therefore induces a morphism of topoi $\varepsilon_S \colon \textsf{Sh}((\textnormal{\textsf{Sm}}/S)_{\textnormal{ét}}) \to \textsf{Sh}((\textnormal{\textsf{Sch}}/S)_{\textnormal{ét}})$ such that the pullback functor $\varepsilon_S^*$ is simply the restriction of sheaves from the big étale site to the lisse-étale site.

    Now let $f \colon T \to S$ be a morphism of schemes. There is a natural pushforward functor $f_* \colon \textsf{Sh}((\textnormal{\textsf{Sm}}/T)_{\textnormal{ét}}) \to \textsf{Sh}((\textnormal{\textsf{Sm}}/S)_{\textnormal{ét}})$ defined by $\Gamma(U,f_* \mathscr{F}) = \Gamma(U \times_S T,\mathscr{F})$. This functor admits a left adjoint $f^*$, but $f^*$ need not be left exact. As a result, the adjunction $(f^*, f_*)$ does not generally define a morphism of topoi.

    However, if we suppose that $f \colon T \to S$ is a \emph{smooth} morphism, then the pullback $f^*$ is given by restriction and is exact. In this case, the adjunction does indeed define a morphism of topoi $f \colon \textsf{Sh}((\textnormal{\textsf{Sm}}/T)_{\textnormal{ét}}) \to \textsf{Sh}((\textnormal{\textsf{Sm}}/S)_{\textnormal{ét}})$.
\end{remark}

\begin{lemma}\label{small fppf cohomology}
    Let $f \colon T \to S$ be a smooth morphism of schemes. Then there is an isomorphism of cohomological $\delta$-functors
    $\normalfont\textrm{H}^i((\textsf{Sch}/S)_\text{ét};T,-)\simeq \textrm{H}^i((\textsf{Sm}/S)_\text{ét};T,\varepsilon_S^*-)$.
\end{lemma}

\begin{proof}
Consider the following commutative diagram, in which all pullback functors are given by restriction:
\[\begin{tikzcd}
	{\textsf{Ab}((\textsf{Sch}/S)_\text{ét})} & {\textsf{Ab}((\textsf{Sm}/S)_\text{ét})} \\
	{\textsf{Ab}((\textsf{Sch}/T)_\text{ét})} & {\textsf{Ab}((\textsf{Sm}/T)_\text{ét}).}
	\arrow["{\varepsilon_S^*}", from=1-1, to=1-2]
	\arrow["{f^*}", from=1-1, to=2-1]
	\arrow["{f^*}", from=1-2, to=2-2]
	\arrow["{\varepsilon_T^*}", from=2-1, to=2-2]
\end{tikzcd}\]
By \cite[Tag \href{https://stacks.math.columbia.edu/tag/03F3}{03F3}]{Stacks}, there is a natural isomorphism
\[\textrm{H}^i((\textsf{Sch}/S)_\text{ét};T,-)\simeq \textrm{H}^i((\textsf{Sch}/T)_\text{ét};T,f^*-).\]
Using a criterion due to Gabber \cite[Thm.\ A.6]{olsson} (see also \cite[Tag \href{https://stacks.math.columbia.edu/tag/0GR0}{0GR0}]{Stacks}), the right-hand side coincides with $\textrm{H}^i((\textsf{Sm}/T)_\text{ét}; T, \varepsilon_T^* f^*-)$. Altogether, we obtain a natural isomorphism
\[\textrm{H}^i((\textsf{Sch}/S)_\text{ét};T,-)\simeq \textrm{H}^i((\textsf{Sm}/T)_\text{ét}; T, f^*\varepsilon_S^*-).\]

It remains to identify $\textrm{H}^i((\textsf{Sm}/T)_\text{ét}; T, f^*-)$ with $\textrm{H}^i((\textsf{Sm}/S)_\text{ét}; T, -)$. This step is subtler than it may appear, as the slice site $(\textsf{Sm}/S)_\text{ét}/T$ is not equivalent to $(\textsf{Sm}/T)_\text{ét}$. The desired natural isomorphism follows from the commutativity of the diagram:
\[\begin{tikzcd}
	{\textsf{Ab}((\textsf{Sm}/S)_\text{ét})} && {\textsf{Ab}((\textsf{Sm}/S)_\text{ét}/T)} && {\textsf{Ab}((\textsf{Sm}/T)_\text{ét}).} \\
	\\
	&& {\textsf{Ab}}
	\arrow[from=1-1, to=1-3]
	\arrow["{f^*}", bend left, from=1-1, to=1-5]
	\arrow["{\textrm{H}^i((\textsf{Sm}/S)_\text{ét};T,-)}"', bend right=25, from=1-1, to=3-3]
	\arrow[from=1-3, to=1-5]
	\arrow["{\textrm{H}^i((\textsf{Sm}/S)_\text{ét}/T;T,-)}", from=1-3, to=3-3]
	\arrow["{\textrm{H}^i((\textsf{Sm}/T)_\text{ét};T,-)}", bend left=25, from=1-5, to=3-3]
\end{tikzcd}\]
The left triangle commutes by another application of \cite[Tag \href{https://stacks.math.columbia.edu/tag/03F3}{03F3}]{Stacks}, whereas the right triangle commutes by Gabber's criterion \cite[Thm.\ A.6]{olsson}.
\end{proof}

\begin{remark}
There are two a priori possible interpretations of Corollary~\ref{small fppf site}. If one follows Definition~\ref{def ext} strictly, then the group $\operatorname{Ext}^i_T(G,\mathscr{A})$ computed on the lisse-étale site $(\mathsf{Sm}/S)_{\text{ét}}$ is, by definition, the value at $\mathscr{A}$ of the $i$-th right derived functor of
\[\operatorname{Hom}_T(G,-)\colon \mathsf{Ab}((\mathsf{Sm}/S)_{\text{ét}}/T)\to \mathsf{Ab}.\]
Alternatively, one may naturally consider the value at $\mathscr{A}$ of the $i$-th right derived functor of
\[\operatorname{Hom}_T(G,-)\colon \mathsf{Ab}((\mathsf{Sm}/T)_{\text{ét}})\to \mathsf{Ab}.\]

Since the slice site $(\mathsf{Sm}/S)_{\text{ét}}/T$ is not, in general, equivalent to $(\mathsf{Sm}/T)_{\text{ét}}$ unless $T$ is étale over $S$, it is not immediate that these two constructions yield the same groups. Nevertheless, the proof of Lemma~\ref{small fppf cohomology} shows that they do agree.
\end{remark}

The following corollary appeared in the first version of this paper under more restrictive hypotheses. We thank T.~Suzuki for explaining how to prove it in the generality stated here.

\begin{corollary}\label{Ext commutes with filtered colimits}
Let $S$ be a quasi-compact and quasi-separated scheme, let $(T_\lambda)_{\lambda\in \Lambda}$ be a cofiltered system of quasi-compact and quasi-separated $S$-schemes with affine transition maps, and let $T=\lim T_\lambda$. Then, for any locally finitely presented commutative group schemes $G$ and $H$ over $S$, the natural map
\[\underset{\lambda\in \Lambda}{\operatorname{colim}}\operatorname{Ext}^i_{T_\lambda}(G,H)\to \operatorname{Ext}^i_T(G,H)\]
is an isomorphism for all $i\geq 0$.
\end{corollary}

\begin{proof}
Since filtered colimits of abelian groups are exact, the Breen--Deligne spectral sequence reduces the claim to the analogous statement for cohomology groups:
\[\underset{\lambda\in \Lambda}{\operatorname{colim}}\: \textrm{H}^i((\textsf{Sch}/T_\lambda)_\text{fppf};G_{T_\lambda},H_{T_\lambda})\to \textrm{H}^i((\textsf{Sch}/T)_\mathrm{fppf};G_T,H_T)\]
is an isomorphism for all $i\geq 0$. 

To analyze this, consider a scheme $X$ and the site $(\textsf{LFP}/X)_\text{fppf}$ of locally finitely presented $X$-schemes with the fppf topology. The inclusion $\textsf{LFP}/X\hookrightarrow \textsf{Sch}/X$ is both continuous and cocontinuous, inducing a morphism of topoi
\[\alpha_X \colon \textsf{Sh}((\textsf{LFP}/X)_\text{fppf}) \to \textsf{Sh}((\textsf{Sch}/X)_\text{fppf})\]
whose pullback functor $\alpha_X^*$ is given by restriction. Gabber's criterion \cite[Thm.\ A.6]{olsson} then shows that the induced map on cohomology
\[\textrm{H}^i((\textsf{Sch}/X)_\text{fppf};X,-)\to\textrm{H}^i((\textsf{LFP}/X)_\text{fppf};X,-)\]
is an isomorphism functorial in $X$. 

Taking $X$ to be $G_T$, and using \cite[Tag \href{https://stacks.math.columbia.edu/tag/03F3}{03F3}]{Stacks}, we obtain a natural isomorphism
\[\textrm{H}^i((\textsf{Sch}/T)_\text{fppf};G_T,H_T)\to\textrm{H}^i((\textsf{LFP}/T)_\text{fppf};G_T,H_T).\]
Thus it suffices to prove that
\[\underset{\lambda\in \Lambda}{\operatorname{colim}}\: \textrm{H}^i((\textsf{LFP}/T_\lambda)_\text{fppf};G_{T_\lambda},H_{T_\lambda})\to \textrm{H}^i((\textsf{LFP}/T)_\text{fppf};G_T,H_T)\]
is an isomorphism for all $i\geq 0$.

Finally, as observed by Grothendieck in the proof of \cite[Lem.\ 11.1]{GrothendieckBrauerIII}, the arguments of \cite[Exp.\ VII, Cor.\ 5.9]{SGA42} apply verbatim in this context, yielding the desired result.
\end{proof}

\section{The Barsotti--Weil formula}
For an abelian $S$-scheme $A$, the abelian fppf sheaf $\underline{\operatorname{Ext}}^1(A,\mathbb{G}_m)$ is represented by the dual abelian scheme of $A$. This result is known as the \emph{Barsotti--Weil formula}. In this section, we present a generalization of it.

For the statement of the next theorem, we introduce the following notation. Given abelian groups  $\mathscr{G}$ and $\mathscr{A}$ in a topos $\normalfont\textsf{X}$, we write $\underline{\operatorname{Mor}}_e(\mathscr{G},\mathscr{A})$ for the sheaf assigning to each object $T$ of $\normalfont\textsf{X}$ the group of morphisms $\mathscr{G}_T \to \mathscr{A}_T$ in $\textsf{X}/T$ that preserve the zero-sections. (The morphisms $\mathscr{G}_T \to \mathscr{A}_T$ are not assumed to preserve the group structure.)

\begin{theorem}[Generalized Barsotti--Weil formula]\label{generalized BW}
    Let $\mathscr{G}$ and $\mathscr{A}$ be an abelian groups in a topos $\normalfont\textsf{X}$ such that $\underline{\operatorname{Mor}}_e(\mathscr{G}^n,\mathscr{A}) = 0$ for $n=1,2,3$. Then, for any object $T$ of $\normalfont\textsf{X}$, the natural maps
    \[\underline{\operatorname{Ext}}^1(\mathscr{G},\mathscr{A})(T)\leftarrow\operatorname{Ext}^1_T(\mathscr{G},\mathscr{A})\to \normalfont\textrm{H}^1_m(\mathscr{G}_T,\mathscr{A}_T)\]
    are isomorphisms. 
    
    In particular, this conclusion applies when $\normalfont\textsf{X}$ is the topos of sheaves on $(\normalfont\textsf{Sch}/S)_\text{fppf}$, $\mathscr{G}$ is represented by an abelian scheme over $S$, and $\mathscr{A}$ is represented either by an affine commutative group scheme over $S$ or by a quasi-coherent $\mathcal{O}_S$-module.
\end{theorem}

\begin{proof}
Since $\underline{\operatorname{Hom}}(\mathscr{G},\mathscr{A})$ is a subsheaf of $\underline{\operatorname{Mor}}_e(\mathscr{G},\mathscr{A})$, Proposition~\ref{sheafification map} implies that, for every object $T$ of $\normalfont\textsf{X}$, the sheafification map
    \[\operatorname{Ext}^1_T(\mathscr{G},\mathscr{A})\to \underline{\operatorname{Ext}}^1(\mathscr{G},\mathscr{A})(T)\]
    is an isomorphism. On the other hand, according to Proposition~\ref{short exact sequence computing extensions}, the natural map
    \[\operatorname{Ext}^1_T(\mathscr{G},\mathscr{A})\to \textrm{H}^1_m(\mathscr{G}_T,\mathscr{A}_T)\]
    is an isomorphism whenever the cohomology groups $\normalfont\textrm{H}^2_s(\mathscr{G}_T,\mathscr{A}_T)$ and $\normalfont\textrm{H}^3_s(\mathscr{G}_T,\mathscr{A}_T)$ both vanish. Since every morphism $\mathscr{G}^n_T \to \mathscr{A}_T$, for $n=2,3$, in $\normalfont\textsf{X}/T$ preserving the zero-sections is trivial, Remark~\ref{normalized BD} shows that these cohomology groups are computed by the complex
    \[0\to 0\to 0\to 0,\]
that is manifestly exact.

We now prove that $\underline{\operatorname{Mor}}_e(\mathscr{G},\mathscr{A})$ vanishes in the aforementioned cases. Let $p\colon A \to S$ be an abelian scheme and $q\colon B \to S$ an affine commutative group scheme. By the universal property of the relative spectrum, giving a morphism of $S$-schemes $A \to B$ is equivalent to specifying a morphism of quasi-coherent $\mathcal{O}_S$-algebras
\[q_{*}\mathcal{O}_{B}\to p_{*}\mathcal{O}_{A}\simeq \mathcal{O}_S,\]
where the final isomorphism follows from \cite[Tag \href{https://stacks.math.columbia.edu/tag/0E0L}{0E0L}]{Stacks}. Thus, any morphism of $S$-schemes $A \to B$ is necessarily constant. In particular, if the morphism preserves the zero-sections, it must be trivial.

Finally, let $M$ be a quasi-coherent sheaf on $S$. By the Yoneda lemma, a morphism of sheaves of sets $\varphi\colon A \to M$ on $(\textsf{Sch}/S)_\mathrm{fppf}$ is determined by a section $s \in \Gamma(A, p^*M)$. Given such a section, the morphism $\varphi$ sends a point $a \in A(T)$ to the pullback $a^*s \in \Gamma(T, a^*p^*M)$. In particular, if $e\colon S \to A$ denotes the zero-section of $A$, then $\varphi$ preserves the zero-sections if and only if $e^*s = 0$ in $\Gamma(S, M)$.

We claim that the pullback $e^*\colon \Gamma(A,p^* M)\to \Gamma(S,M)$ is an isomorphism, from which the result follows. Since $p$ is flat, the projection formula \cite[Tag \href{https://stacks.math.columbia.edu/tag/08EU}{08EU}]{Stacks} gives an isomorphism
\[\mathsf{R}p_* p^* M\simeq \mathsf{R}p_*(\mathcal{O}_A\otimes^\mathsf{L}_{\mathcal{O}_A}p^* M)\simeq \mathsf{R}p_*\mathcal{O}_A\otimes^\mathsf{L}_{\mathcal{O}_S} M.\]
As a result, the Künneth spectral sequence, whose second page is
\[E_2^{n,m}=\bigoplus_{i+j=m}\underline{\operatorname{Tor}}^n_{\mathcal{O}_S}(\mathsf{R}^ip_*\mathcal{O}_A,\mathscr{H}^j(M)),\]
converges to $\mathsf{R}^{n+m}p_* p^* M$. Now, by \cite[Thm.\ 27.203]{GWII}, each $\mathsf{R}^i p_* \mathcal{O}_A$ is a finite locally free $\mathcal{O}_S$-module. Consequently, these $\underline{\operatorname{Tor}}^n_{\mathcal{O}_S}$ vanish for all $n > 0$, and the spectral sequence degenerates. We thus obtain an isomorphism $p_* p^* M\simeq M$, and the claim follows by taking global sections.
\end{proof}

Although the Breen--Deligne resolution simplifies the proof of the preceding theorem, and consequently those of the corollaries in this section, it is by no means essential. We therefore also provide a direct proof that the natural map $\operatorname{Ext}^1_T(\mathscr{G},\mathscr{A})\to \textrm{H}^1_m(\mathscr{G}_T,\mathscr{A}_T)$ is an isomorphism.

\begin{proof}[Alternative proof of Theorem~\ref{generalized BW}]
As usual, the group $\operatorname{Ext}^1_T(\mathscr{G},\mathscr{A})$ classifies extensions of $\mathscr{G}_T$ by $\mathscr{A}_T$ in $\textsf{Ab}(\textsf{X}/T)$. Consider an element of $\operatorname{Ext}^1_T(\mathscr{G},\mathscr{A})$ represented by an exact sequence
\[0\to \mathscr{A}_T\to \mathscr{E}\xrightarrow{\ \pi\ } \mathscr{G}_T\to 0.\]
Such an extension lies in the kernel of the map $\operatorname{Ext}^1_T(\mathscr{G},\mathscr{A})\to \textrm{H}^1_m(\mathscr{G}_T,\mathscr{A}_T)$ if and only if the associated $\mathscr{A}_T$-torsor over $\mathscr{G}_T$ is trivial. In concrete terms, this means there exists a morphism $\sigma\colon \mathscr{G}_T \to \mathscr{E}$ in $\textsf{X}/T$ satisfying $\pi \circ \sigma = \operatorname{id}_{\mathscr{G}_T}$. To prove that the map $\operatorname{Ext}^1_T(\mathscr{G},\mathscr{A})\to \textrm{H}^1_m(\mathscr{G}_T,\mathscr{A}_T)$ is injective, it suffices to show that there exists such a section $\sigma$ that is a morphism of groups. After composing $\sigma$ with translation by a global section of $\mathscr{A}_T$, we may assume that $\sigma$ preserves the zero-sections.

Consider the morphism $\theta\colon \mathscr{G}^2_T\to \mathscr{E}$ given by $(x,y)\mapsto \sigma(xy)\sigma(y)^{-1}\sigma(x)^{-1}$. We claim that this map factors through $\mathscr{A}_T$. Indeed, an element of $\mathscr{E}$ lies in $\mathscr{A}_T$ if and only if its image by $\pi$ vanishes and we have that
\begin{align*}
   \pi(\theta(x,y))&=\pi(\sigma(xy)\sigma(y)^{-1}\sigma(x)^{-1}) \\
   &= \pi(\sigma(xy))\pi(\sigma(y))^{-1}\pi(\sigma(x))^{-1}\\
   &= xyy^{-1}x^{-1}=e.
\end{align*}
By assumption, every morphism $\mathscr{G}_T^2\to \mathscr{A}_T$ in $\textsf{X}/T$ preserving the zero-sections must be trivial. This implies that $\theta$ is the constant map to the zero-section of $\mathscr{A}_T$, and hence that $\sigma$ is a morphism of groups.

To prove that $\operatorname{Ext}^1_T(\mathscr{G}, \mathscr{A}) \to \textrm{H}^1_m(\mathscr{G}_T, \mathscr{A}_T)$ is surjective, let $\pi\colon P\to \mathscr{G}_T$ be a multiplicative $\mathscr{A}_T$-torsor. The fact that the torsor is multiplicative exactly says that there is a map $m_P\colon P \times_T P \to P$ lying over the group law on $\mathscr{G}_T$ and equivariant with respect to the $\mathscr{A}_T$-action on both factors:
\[m_P(a\cdot x, y) = m_P(x, a\cdot y) = a\cdot m_P(x, y).\]
We will modify this map and show that it then yields a group structure on $P$, making it into an extension of $\mathscr{G}_T$ by $\mathscr{A}_T$.

First, we note that restricting the multiplicativity condition along the zero-section of $\mathscr{G}_T^2$ shows that the resulting $\mathscr{A}$-torsor over $T$ is trivial. That is, there is a point $e_P \in P(T)$ lying above the zero-section of $\mathscr{G}_T$. Consider the map $R\colon P \to P$, defined as $R(x)=m_P(x, e_P)$. This morphism lies above the identity map of $\mathscr{G}_T$, hence we have $R(x) = \rho(x)\cdot x$ for some map $\rho\colon P \to \mathscr{A}_T$. The $\mathscr{A}_T$-equivariance of $m_P$ implies that $\rho$ is $\mathscr{A}_T$-invariant, and therefore descends to a morphism $\mathscr{G}_T \to \mathscr{A}_T$, which we also denote by $\rho$.

If we replace $m_P$ by the map $(x, y) \mapsto (-\rho(x))\cdot m_P(x, y)$, then $e_P$ becomes a right identity for $m_P$, and this change does not affect the $\mathscr{A}_T$-equivariance of $m_P$ or the fact that $m_P$ lies above multiplication on $\mathscr{G}_T$. We may therefore assume that $e_P$ is a right identity. Note that it is the unique right identity, because any right identity $e^\prime_P$ must lie above the zero-section of $\mathscr{G}_T$, hence $e^\prime_P = a\cdot e_P$ for some $a \in \mathscr{A}_T$, so one has
\[e_P = m_P(e_P, e^\prime_P) = m_P(e_P, a\cdot e_P) = a\cdot m_P(e_P, e_P) = a\cdot e_P,\]
hence $a$ is the zero-section of $\mathscr{A}_T$ and $e'_P = e_P$.

Next we verify the existence of left inverses. For this, we use the multiplicativity again. The isomorphism $m^*P \simeq \operatorname{pr}_1^*P \wedge\operatorname{pr}_2^*P$ pulled back along the map 
\[(\operatorname{id},-\operatorname{id})\colon \mathscr{G}_T \to \mathscr{G}_T^2,\]
together with the already-used fact that the restriction of $P$ to $T$ is trivial (also a consequence of multiplicativity), shows that $[-1]^*P \simeq  -P$ as $\mathscr{A}_T$-torsors over $\mathscr{G}_T$. This implies the existence of a map $\operatorname{inv}_P\colon P \to P$ lying above inversion on $\mathscr{G}_T$ and satisfying 
\begin{equation} \label{inv}
\operatorname{inv}_P((-b)\cdot x) = b\cdot \operatorname{inv}_P(x) \tag{$\ast$}
\end{equation}
for $b \in \mathscr{A}_T$.
Note that $\operatorname{inv}_P(e_P) = a\cdot e_P$ for some $a \in \mathscr{A}_T$. If we replace $\operatorname{inv}_P$ by the map $(-a)\cdot \operatorname{inv}_P$, then this new map still lies above inversion on $\mathscr{G}_T$ and satisfies \eqref{inv}, but now additionally satisfies $\operatorname{inv}_P(e_P) = e_P$.

As above, the obstacle to the identity $m_P(\operatorname{inv}_P(x), x) = e_P$ holding is a map $\iota\colon \mathscr{G}_T \to \mathscr{A}_T$ satisfying $m_P(\operatorname{inv}_P(x), x) = \iota(\pi(x))\cdot e_P$. Note that the equation $\operatorname{inv}_P(e_P) = e_P$ implies that $\iota$ preserves zero-sections. We replace $m_P$ by the map $(x,y)\mapsto (-\iota(\pi(y)))\cdot m_p(x, y)$.\footnote{Under the assumption that every map $\mathscr{G}_T\to \mathscr{A}_T$ is constant, this step is unecessary. We do it to make clear that the map $\operatorname{Ext}^1_T(\mathscr{G},\mathscr{A})\to \textrm{H}^1_m(\mathscr{G}_T,\mathscr{A}_T)$ is an isomorphism if every map $\mathscr{G}_T^n\to \mathscr{A}_T$, for $n=2,3$, is constant. This is also clear from the proof based on the Breen--Deligne resolution.} This preserves the property that $m_P$ is $\mathscr{A}_T$-equivariant, lies above multiplication on $\mathscr{G}_T$, and has $e_P$ has a right identity, but the new $m_P$ now has $\operatorname{inv}_P$ as a left inverse.

Next we consider commutativity of $m_P$. Because $\mathscr{G}_T$ is commutative, we have identically
\[m_P(x, y) = \gamma(x, y)\cdot m_P(y, x),\]
for some morphism $\gamma\colon P^2 \to \mathscr{A}_T$. The $\mathscr{A}_T$-equivariance of $m_P$ implies that $\gamma$ is $\mathscr{A}_T$\nobreakdash-invariant, hence descends to a map (which, as before, we still denote by $\gamma$) $\mathscr{G}_T^2 \to \mathscr{A}_T$. The only such morphism is constant, and the fact that $e_P$ is a right identity for $m_P$ implies that $\gamma$ preserves zero-sections, so $\gamma$ is trivial and $m_P$ is commutative. An exactly analogous argument, using the triviality of pointed maps $\mathscr{G}_T^3 \rightarrow \mathscr{A}_T$, shows that $m_P$ is associative.

Thus $m_P$ defines a group law on $P$, and the map $\pi\colon P \to \mathscr{G}_T$ is a homomorphism. The kernel is the image of the map $\mathscr{A}_T \to P$, $a \mapsto a\cdot e_P$. This is an injective map, and it is a homomorphism due to the $\mathscr{A}_T$-equivariance of $m_P$. Finally, the $\mathscr{A}_T$-equivariance of $m_P$ also shows that $P$, considered as a $\mathscr{A}_T$-torsor over $\mathscr{G}_T$ via translation by $\mathscr{A}_T$, is in fact just the torsor with which we began.
\end{proof}

\begin{remark}
The maps $\mathscr{G}_T^2, \mathscr{G}_T^3 \to \mathscr{A}_T$ appearing in the above proof yield explicit descriptions of the morphism
\[\textrm{H}^1_m(\mathscr{G}_T, \mathscr{A}_T) \to \textrm{H}^3_s(\mathscr{G}_T, \mathscr{A}_T)\]
appearing in Proposition~\ref{short exact sequence computing extensions}, which obstructs multiplicative torsors from arising via extensions. Recall that $\textrm{H}^3_s(\mathscr{G}_T, \mathscr{A}_T)$ is defined as a subquotient of $\Gamma(\mathscr{G}_T^3, \mathscr{A}_T) \oplus \Gamma(\mathscr{G}_T^2, \mathscr{A}_T)$. The first summand measures the failure of associativity, while the second measures the failure of commutativity for the map $m_P$ described in the above proof.
\end{remark}

\begin{remark}\label{noncommutative remark}
The assumption that $\mathscr{G}_T$ is commutative is not essential: even without it, the same arguments show---under the hypothesis that there are no nonconstant morphisms $\mathscr{G}_T^n \to \mathscr{A}_T$ for $n = 2, 3$---that the map from central Yoneda extensions of $\mathscr{G}_T$ by $\mathscr{A}_T$ to $\textrm{H}^1_m(\mathscr{G}_T, \mathscr{A}_T)$ is a bijection. See \cite[Exp.\ VII, Prop.\ 1.3.5]{raynaud2006groupes} and \cite[Thm.\ 3.3]{bruin} for noncommutative analogues of Proposition~\ref{short exact sequence computing extensions}, proved using similar arguments.
\end{remark}

Before explaining some consequences of this generalized Barsotti--Weil formula, we prove a Künneth formula for the cohomology of an abelian scheme over an arbitrary ring.

\begin{lemma}\label{Kunneth lemma}
    Let $A$ be an abelian scheme over a ring $R$. Then the $R$-module $\normalfont\textrm{H}^n(A^2,\mathcal{O}_{A^2})$ is the direct sum of $\normalfont\textrm{H}^{i}(A,\mathcal{O}_{A})\otimes_R\textrm{H}^{j}(A,\mathcal{O}_{A})$ for $i+j=n$. In particular, $\normalfont\textrm{H}^1(A^2,\mathcal{O}_{A^2})$ is isomorphic to $\normalfont\textrm{H}^{1}(A,\mathcal{O}_{A})\oplus \textrm{H}^{1}(A,\mathcal{O}_{A})$.
\end{lemma}

\begin{proof}
    Since $A$ is flat over $R$, we have $\mathsf{R}\Gamma(A^2,\mathcal{O}_{A^2})\simeq \mathsf{R}\Gamma(A,\mathcal{O}_A)\otimes^\mathsf{L}_R \mathsf{R}\Gamma(A,\mathcal{O}_A)$, as shown in \cite[Tag \href{https://stacks.math.columbia.edu/tag/0FLQ}{0FLQ}]{Stacks}. Akin to the first proof of Theorem~\ref{generalized BW}, the associated Künneth spectral sequence degenerates, yielding the desired isomorphism.
\end{proof}

\begin{corollary}\label{unipotent BW}
	Let $A$ be an abelian scheme over $S$, and let $M$ be a quasi-coherent $\mathcal{O}_S$-module. For a morphism of schemes $f\colon T\to S$, the natural maps
	\[\underline{\operatorname{Ext}}^1(A,M)(T)\leftarrow\operatorname{Ext}^1_T(A,M)\to \normalfont\textrm{H}^1(A_T,f^*M)/p_T^*{\textrm{H}}^1(T, f^*M),\]
     where $p\colon A \to S$ is the structure morphism and $p_T\colon A_T \to T$ its base change, are isomorphisms. Moreover, the sheaf $\underline{\operatorname{Ext}}^1(A,M)$ is isomorphic to $\operatorname{Lie}(A^\vee) \otimes_{\mathcal{O}_S} M$, where $\operatorname{Lie}(A^\vee)$ denotes the Lie algebra of the dual abelian scheme.
\end{corollary}

\begin{proof}
    We first claim that there is an isomorphism
    \[\Gamma(T,f^*(\operatorname{Lie}(A^\vee) \otimes_{\mathcal{O}_S} M))\simeq \textrm{H}^1(A_T,f^*M)/p_T^*{\textrm{H}}^1(T, f^* M),\]
    that is functorial in $T$. By \cite[Prop.\ 27.122]{GWII}, we have an isomorphism of quasi-coherent sheaves $\operatorname{Lie}(A^\vee)\simeq \mathsf{R}^1p_* \mathcal{O}_A$ and, by \cite[Thm.\ 27.203]{GWII}, the formation of the latter is compatible with base change. Hence,
    \[\Gamma(T,f^*(\operatorname{Lie}(A^\vee) \otimes_{\mathcal{O}_S} M))\simeq \Gamma(T,f^* \mathsf{R}^1p_* \mathcal{O}_A\otimes_{\mathcal{O}_T} f^*M)\simeq \Gamma(T,\mathsf{R}^1p_{T,*} \mathcal{O}_{A_T}\otimes_{\mathcal{O}_T} f^*M).\]
    As in the proof of Theorem~\ref{generalized BW}, this is the same as $\Gamma(T,\mathsf{R}^1p_{T,*} p_T^* f^*M)$. Using the Leray spectral sequence, we obtain an exact sequence
\[\begin{tikzcd}
	0 & {\textrm{H}^1(T,f^*M)} & {\textrm{H}^1(A_T,f^*M)} & {\Gamma(T,\mathsf{R}^1p_{T,*}p_T^* f^*M)}\ar[draw=none]{dll}[name=X, anchor=center]{}\ar[rounded corners,
            to path={ -- ([xshift=2ex]\tikztostart.east)
                      |- (X.center) \tikztonodes
                      -| ([xshift=-2ex]\tikztotarget.west)
                      -- (\tikztotarget)}]{dll}[at end]{} \\
	& {\textrm{H}^2(T,f^*M)} & {\textrm{H}^2(A_T, f^*M).}
	\arrow[from=1-1, to=1-2]
	\arrow[from=1-2, to=1-3]
	\arrow[from=1-3, to=1-4]
	\arrow[from=2-2, to=2-3]
\end{tikzcd}\]
    Since the structure map $p_T\colon A_T\to T$ has a section, the induced pullback morphism $\textrm{H}^2(T,f^*M)\to \textrm{H}^2(A_T, f^*M)$ is injective, thus establishing the desired isomorphism.

	The remaining content of this corollary is that the natural map
    \[\textrm{H}^1_m(A_T, f^*M)\to \textrm{H}^1(A_T, f^*M)/p_T^*{\textrm{H}}^1(T, f^*M)\]
    is an isomorphism. By Theorem~\ref{generalized BW} and the discussion above, both sides are fppf sheaves on $T$ and the morphism above is functorial. In particular, we may assume that $T$ is affine.
    
    Let $m\colon A_T\times_T A_T\to A_T$ denote the group law, and let $\operatorname{pr}_1,\operatorname{pr}_2\colon A_T\times_T A_T\to A_T$ be the natural projections. From the isomorphism $\textrm{H}^i(A_T,f^*M)\simeq \textrm{H}^i(A_T,\mathcal{O}_{A_T})\otimes_{\mathcal{O}_T}f^*M$ seen above (and its analogue for $A_T^2$), together with the Küneth formula, we obtain that the map
	\begin{align*}
	\textrm{H}^1(A_T,f^*M)\oplus \textrm{H}^1(A_T,f^*M)&\to\textrm{H}^1(A_T^2,f^*M)\\
	(x,y) &\mapsto \operatorname{pr}_1^*x+\operatorname{pr}_2^*y
	\end{align*}
	is an isomorphism. In particular, given an element $z$ of $\textrm{H}^1(A_T,f^*M)$, there exist $x,y\in \textrm{H}^1(A_T,f^*M)$ satisfying
	\[\operatorname{pr}_1^*x+\operatorname{pr}_2^*y=m^*z.\]
    Restricting this along $\operatorname{id}\times \, e\colon A_T\simeq  A_T\times_T T\to A_T\times_T A_T$, we obtain that $x-z$ lies in $p_T^*{\textrm{H}}^1(T, f^*M)$, which vanishes due to the assumption that $T$ is affine. Similarly, restricting it along $e\times \operatorname{id}$ we obtain that $y=z$; proving that $z$ is multiplicative.
\end{proof}

\begin{corollary}\label{usual BW}
    Let $A$ be an abelian scheme over $S$. Then the abelian fppf sheaf $\underline{\operatorname{Ext}}^1(A,\mathbb{G}_m)$ is representable by the dual abelian scheme $A^\vee$.
\end{corollary}

\begin{proof}
Let $p\colon A \to S$ be the structure map, and denote its base change to an $S$-scheme $T$ by $p_T\colon A_T \to T$. Let $m\colon A_T \times_T A_T \to A_T$ be the group law, and let $\operatorname{pr}_1, \operatorname{pr}_2 \colon A_T \times_T A_T \to A_T$ denote the natural projections. According to \cite[Prop.~27.161]{GWII}, the group $A^\vee(T)$ is isomorphic to the kernel of the map
\[
\operatorname{pr}_1^* + \operatorname{pr}_2^*-m^* \colon \frac{\operatorname{Pic}(A_T)}{p_T^* \operatorname{Pic}(T)} \to \frac{\operatorname{Pic}(A_T^2)}{(p_T \times p_T)^* \operatorname{Pic}(T)}.
\]

Note that, for each of the morphisms $m$, $\operatorname{pr}_1$, and $\operatorname{pr}_2$, the composition $A_T^2 \to A_T \xrightarrow{p_T} T$ coincides with $p_T \times p_T$. This implies that the left square in the following diagram commutes:
\[\begin{tikzcd}
	0 & {\operatorname{Pic}(T)} & {\operatorname{Pic}(A_T)} & {\operatorname{Pic}(A_T)/p_T^*\operatorname{Pic}(T)} & 0 \\
	0 & {\operatorname{Pic}(T)} & {\operatorname{Pic}(A_T^2)} & {\operatorname{Pic}(A^2_T)/(p_T\times p_T)^*\operatorname{Pic}(T)} & 0.
	\arrow[from=1-1, to=1-2]
	\arrow[from=1-2, to=1-3]
	\arrow[equals, from=1-2, to=2-2]
	\arrow[from=1-3, to=1-4]
	\arrow["{\operatorname{pr}_1^*+\operatorname{pr}_2^*-m^*}", from=1-3, to=2-3]
	\arrow[from=1-4, to=1-5]
	\arrow[from=1-4, to=2-4]
	\arrow[from=2-1, to=2-2]
	\arrow[from=2-2, to=2-3]
	\arrow[from=2-3, to=2-4]
	\arrow[from=2-4, to=2-5]
\end{tikzcd}\]
The snake lemma then yields a functorial isomorphism $\textrm{H}^1_m(A_T,\mathbb{G}_m)\simeq A^\vee(T)$, and the desired result follows from Theorem~\ref{generalized BW}.
\end{proof}

\section{Higher extensions of abelian schemes}\label{higher exts}
After computing the extension sheaves $\underline{\operatorname{Ext}}^1(A,\mathbb{G}_a) \simeq \operatorname{Lie}(A^\vee)$ and $\underline{\operatorname{Ext}}^1(A,\mathbb{G}_m) \simeq A^\vee$, we now study their higher analogs $\underline{\operatorname{Ext}}^2(A,\mathbb{G}_a)$ and $\underline{\operatorname{Ext}}^2(A,\mathbb{G}_m)$. Since proving the vanishing of the former is both simpler and a necessary step toward proving the vanishing of the latter, we begin with it.

\begin{theorem}\label{extension by quasi-coherent sheaf}
    Let $A$ be an abelian scheme over a ring $R$, and let $M$ be an $R$-module. Then $\operatorname{Ext}^2_R(A,M)=0$.
\end{theorem}

\begin{proof}
Consider the structure map $p\colon A\to \operatorname{Spec}R$. As seen in the proof of Corollary~\ref{unipotent BW}, the cohomology group $\textrm{H}^j(A,p^*M)$ is isomorphic to $\textrm{H}^j(A,\mathcal{O}_A)\otimes_R M$ for all $j$. Next, we compute $\operatorname{Ext}^2_R(A,M)$ using the Breen--Deligne spectral sequence, as in Proposition~\ref{short exact sequence computing extensions}. To simplify notation, let $\textrm{H}^{i,j}$ denote the $R$-module $\textrm{H}^i(A^j,\mathcal{O}_{A^j})\otimes_R M$. We explicitly write the beginning of the first page.
\[\begin{tikzcd}[ampersand replacement=\&,column sep=small,row sep=tiny]
	{\textrm{H}^{2,1}} \& {\textrm{H}^{2,2}} \& {\textrm{H}^{2,3}\oplus \textrm{H}^{2,2}} \& {\textrm{H}^{2,4}\oplus \textrm{H}^{2,3}\oplus \textrm{H}^{2,3}\oplus \textrm{H}^{2,2}\oplus \textrm{H}^{2,1}} \\
	{\textrm{H}^{1,1}} \& {\textrm{H}^{1,2}} \& {\textrm{H}^{1,3}\oplus \textrm{H}^{1,2}} \& {\textrm{H}^{1,4}\oplus \textrm{H}^{1,3}\oplus \textrm{H}^{1,3}\oplus \textrm{H}^{1,2}\oplus \textrm{H}^{1,1}} \\
	{\textrm{H}^{0,1}} \& {\textrm{H}^{0,2}} \& {\textrm{H}^{0,3}\oplus \textrm{H}^{0,2}} \& {\textrm{H}^{0,4}\oplus \textrm{H}^{0,3}\oplus \textrm{H}^{0,3}\oplus \textrm{H}^{0,2}\oplus \textrm{H}^{0,1}}
	\arrow[from=1-1, to=1-2]
	\arrow[from=1-2, to=1-3]
	\arrow[from=1-3, to=1-4]
	\arrow[from=2-3, to=2-4]
	\arrow[from=3-3, to=3-4]
	\arrow[from=2-1, to=2-2]
	\arrow[from=2-2, to=2-3]
	\arrow[from=3-1, to=3-2]
	\arrow[from=3-2, to=3-3]
\end{tikzcd}\]
Here, the maps on the left are given by $(m^*-\operatorname{pr}_1^*-\operatorname{pr}_2^*)\otimes \operatorname{id}$, where $m\colon A\times A\to A$ is the group operation. The maps on the middle are given by
\[(\operatorname{pr}_{1,2}^*-(\operatorname{id}\times m)^*+(m\times\operatorname{id})^*-\operatorname{pr}_{2,3}^*,\operatorname{id}^*-\tau^*)\otimes \operatorname{id},\]
where $\tau\colon A\times A\to A\times A$ permutes the factors. Similarly, the maps on the right can be computed using the formulas described just after Proposition~\ref{BD resolution}.

First we check that the bottom row is exact at $\textrm{H}^{0,3}\oplus \textrm{H}^{0,2}$, for which it suffices to check that the map $$\textrm{H}^{0,3}\oplus \textrm{H}^{0,2} \to \textrm{H}^{0,4} \oplus \textrm{H}^{0,1}$$ is injective. Since $\textrm{H}^0(A^j,\mathcal{O}_{A^j})\simeq R$ for all $j$, this map is simply
\[\begin{tikzcd}[ampersand replacement=\&,row sep=0.25em]
	{M^{\oplus 2}} \& {M^{\oplus 2}} \\
	{(x,y)} \& {(-x,y),}
	\arrow[from=1-1, to=1-2]
	\arrow[maps to, from=2-1, to=2-2]
\end{tikzcd}\]
which is clearly injective.

To show that the cohomology of the middle row at $\textrm{H}^{1,2}$ vanishes, it suffices to prove that the map $\textrm{H}^{1,2}\to \textrm{H}^{1,3}$ is injective. According to Lemma~\ref{Kunneth lemma}, we have isomorphisms $\textrm{H}^{1,2}\simeq (\textrm{H}^{1,1})^{\oplus 2}$ and $\textrm{H}^{1,3}\simeq (\textrm{H}^{1,1})^{\oplus 3}$. Now, consider the composition
\[\textrm{H}^{1,1}\to (\textrm{H}^{1,1})^{\oplus 2}\simeq \textrm{H}^{1,2}\to \textrm{H}^{1,3}\simeq (\textrm{H}^{1,1})^{\oplus 3}\to \textrm{H}^{1,1},\]
where the middle map is one of $\{\operatorname{pr}_{1,2}^*,(\operatorname{id}\times m)^*,(m\times\operatorname{id})^*,\operatorname{pr}_{2,3}^*\}$, and the other maps are some choice of the natural inclusions or projections. This composition can be geometrically described as the pullback by
\[A \to A\times_R A\times_R A \to A\times_R A\to A,\]
where the first map is the closed immersion of a factor by restricting to zero on the other factors, the middle map is one of $\{\operatorname{pr}_{1,2},\operatorname{id}\times m,m\times\operatorname{id},\operatorname{pr}_{2,3}\}$, and the last map is a projection.

By considering all possible inclusions $\textrm{H}^{1,1}\to (\textrm{H}^{1,1})^{\oplus 2}$, middle maps, and projections $(\textrm{H}^{1,1})^{\oplus 3}\to \textrm{H}^{1,1}$, we see that the map $\textrm{H}^{1,2}\to \textrm{H}^{1,3}$ acts as
\[\begin{tikzcd}[ampersand replacement=\&,row sep=0.25em]
	{(\textrm{H}^{1,1})^{\oplus 2}} \& {(\textrm{H}^{1,1})^{\oplus 3}} \\
	{(x,y)} \& {(x,0,-y),}
	\arrow[from=1-1, to=1-2]
	\arrow[maps to, from=2-1, to=2-2]
\end{tikzcd}\]
which is clearly injective.

To prove that the morphism $\textrm{H}^{2,1} \to \textrm{H}^{2,2}$ is injective, we could carry out a direct calculation as above. However, we prefer to leverage the computation already performed in \cite[Lem.\ 27.209]{GWII}, which assumes that $R$ is a field and $M=R$. First, we may assume that $R$ is a local ring, since injectivity can be checked after localization at prime ideals. By \cite[Cor.\ 27.204]{GWII}, the cohomology algebra $\textrm{H}^*(A, \mathcal{O}_A)$ is then identified with the exterior algebra on the \emph{free} $R$-module $\textrm{H}^1(A, \mathcal{O}_A)$ of rank $g \colonequals \dim_R(A)$.

Arbitrarily identifying $\textrm{H}^1(A, \mathcal{O}_A)$ with $R^g$, the maps $m^*$ and $\operatorname{pr}_i^*$ on cohomology algebras are then obtained by universal formulas (depending only on $g$) that are completely independent of $R$ or $A$. We therefore have a map
\[f\colon \mathbb{Z}^{\oplus\binom{g}{2}} \to \mathbb{Z}^{\oplus \binom{2g}{2}},\]
and we would like to show that $f \otimes N$ is injective for any $\mathbb{Z}$-module $N$. We know that this holds when $N$ is a field, due to the result \cite[Lem.\ 27.209]{GWII} cited earlier.

To show that $f$ itself is injective, it is enough to verify this after tensoring to $\mathbb{F}_p$ for every (in fact, because $f$ is a map between free finite rank modules, even just a single) prime $p$, and in this case injectivity holds because $\mathbb{F}_p$ is a field. Now we show that $f$ remains injective after tensoring with any $\mathbb{Z}$-module. Letting $C \colonequals \mathrm{coker}(f)$, it suffices (is equivalent, in fact) to show that $C$ is $\mathbb{Z}$-flat, for which it is enough to show that $\operatorname{Tor}^1(C,\mathbb{F}_p) = 0$ for every prime $p$. Because $f$ has flat codomain, this follows from the injectivity of $f \otimes \mathbb{F}_p$. This finishes the proof.
\end{proof}

\begin{corollary}\label{vanishing of extensions by additive group}
    Let $A$ be an abelian scheme over a base scheme $S$, and let $M$ be a quasi-coherent $\mathcal{O}_S$-module. For a morphism of schemes $f\colon T\to S$, there is an isomorphism
    \[\operatorname{Ext}_T^2(A,M)\simeq \normalfont\textrm{H}^1(T,f^*(\operatorname{Lie}(A^\vee)\otimes_{\mathcal{O}_S}M))\]
    that is functorial in $T$. In particular, $\underline{\operatorname{Ext}}^2(A,M)$ vanishes as a sheaf in the big Zariski site of $S$.
\end{corollary}

\begin{proof}
    The sheaf $\underline{\operatorname{Ext}}^2(A, M)$ is the fppf sheafification of the presheaf $T \mapsto \operatorname{Ext}^2_T(A, M)$. Since this presheaf vanishes on affine schemes, it follows that $\underline{\operatorname{Ext}}^2(A, M)$ is zero as a sheaf on the big fppf site of $S$. The desired isomorphism then follows from Proposition~\ref{sheafification map}, which in turn implies that $\underline{\operatorname{Ext}}^2(A, M)$ vanishes as a sheaf on the big Zariski site of $S$.
\end{proof}

We now state the main result of this section.

\begin{theorem}\label{higher vanishing for abelian schemes}
Let $A$ be an abelian scheme over a base scheme $S$, and let $T$ be a $S$-scheme. There is an isomorphism
\[\operatorname{Ext}_T^2(A,\mathbb{G}_m)\simeq \normalfont\textrm{H}^1(T,A^\vee_T)\]
that is functorial in $T$. In particular, $\underline{\operatorname{Ext}}^2(A,\mathbb{G}_m)$ vanishes as a sheaf in the big étale site of $S$.
\end{theorem}

Contrary to the case of the additive group, the group $\operatorname{Ext}_T^2(A,\mathbb{G}_m)$ may fail to vanish even when $T$ is affine. Numerous examples of nontrivial elements in $\textrm{H}^1(T, A^\vee_T)$ can be found in \cite[Chap.\ XIII]{raynaud}. Furthermore, Proposition~\ref{sheafification map} shows that Theorem~\ref{higher vanishing for abelian schemes} actually follows from the weaker statement that $\underline{\operatorname{Ext}}^2(A,\mathbb{G}_m)$ vanishes as a sheaf on the big fppf site of $S$.

The proof of vanishing of $\underline{\operatorname{Ext}}^2(A,\mathbb{G}_m)$ in the fppf topology requires several lemmas. Readers are encouraged to first read the proof in full and refer back to the lemmas as needed.

\begin{lemma}\label{Ext is torsion-free}
Let $A$ an abelian scheme over a base scheme $S$. For every non-zero integer $n$, the multiplication-by-$n$ map on $\underline{\operatorname{Ext}}^2(A,\mathbb{G}_m)$ is injective.
\end{lemma}

\begin{proof}
According to \cite[Prop.\ 27.186]{GWII}, we have a short exact sequence of abelian sheaves
\[0\to A[n]\to A\xrightarrow{\ n\ }A\to 0,\]
where $A[n]$ is a finite locally free group scheme over $S$. Passing to the long exact sequence in cohomology, we obtain the exact sequence
\[\underline{\operatorname{Ext}}^1(A[n],\mathbb{G}_m)\to \underline{\operatorname{Ext}}^2(A,\mathbb{G}_m)\xrightarrow{\ n\ } \underline{\operatorname{Ext}}^2(A,\mathbb{G}_m)\]
and $\underline{\operatorname{Ext}}^1(A[n],\mathbb{G}_m)$ vanishes due to \cite[Exp.\ VIII, Prop.\ 3.3.1]{raynaud2006groupes}.
\end{proof}

\begin{lemma}\label{Ext and reduced schemes}
Let $A$ be an abelian scheme over a noetherian ring $R$. Then the restriction map
\[\operatorname{Ext}^2_R(A,\mathbb{G}_m)\to \operatorname{Ext}^2_{R_\mathrm{red}}(A,\mathbb{G}_m)\]
is injective.
\end{lemma}

\begin{proof}
    Let $I \subset R$ be an arbitrary nilpotent ideal, and let $n$ denote its index of nilpotency. We will show that the map $\operatorname{Ext}^2_R(A,\mathbb{G}_m)\to \operatorname{Ext}^2_{R/I}(A,\mathbb{G}_m)$ is injective, which suffices due to the noetherianity of $R$. Note that the quotient map $R\to R/I$ factors as
    \[R=R/I^n\to R/I^{n-1}\to \dots \to R/I^2\to R/I,\]
    where every map is surjective and has a square-zero kernel. Consequently, we may assume that $n=2$. Moreover, since $\mathbb{G}_m$ is smooth, we may compute extension groups as derived functors on the big étale site.

    Let $i\colon \operatorname{Spec}R/I\to\operatorname{Spec}R$ be the associated closed immersion. The unit of adjunction $\mathbb{G}_{m,R}\to i_*\mathbb{G}_{m,R/I}$ is an epimorphism of abelian sheaves in $(\textsf{Sch}/\operatorname{Spec}R)_\text{ét}$, and we denote by $K$ its kernel. Consequently, we obtain an exact sequence
\[\operatorname{Ext}^2_R(A,K)\to \operatorname{Ext}^2_R(A,\mathbb{G}_m)\to \operatorname{Ext}^2_{R}(A,i_*\mathbb{G}_m).\]
   Since $i_*$ is exact by \cite[Tag \href{https://stacks.math.columbia.edu/tag/04C4}{04C4}]{Stacks}, the groups $\operatorname{Ext}^2_{R}(A,i_*\mathbb{G}_m)$ and $\operatorname{Ext}^2_{R/I}(A,\mathbb{G}_m)$ are naturally isomorphic. As a result, it suffices to prove that $\operatorname{Ext}^2_R(A,K)$ vanishes.

   We now claim that the restrictions of $K$ and $I$ to the lisse-étale site $(\textsf{Sm}/\operatorname{Spec} R)_{\text{ét}}$ are isomorphic. To see this, note that for an $R$-algebra $B$, the sections of $K$ over $\operatorname{Spec} B$ are given by $1 + IB$, which is isomorphic to $IB$ because $I$ is square-zero. On the other hand, the sections of the quasi-coherent sheaf $I$ over $\operatorname{Spec} B$ are $I \otimes_R B$. There is thus a natural morphism $I \to K$, which becomes an isomorphism when restricted to flat $R$-schemes. By Corollary~\ref{small fppf site}, this induces an isomorphism $\operatorname{Ext}^2_R(A, I) \simeq \operatorname{Ext}^2_R(A, K)$, and the former group vanishes by Theorem~\ref{extension by quasi-coherent sheaf}.
\end{proof}

We conclude by stating the final lemma, whose proof will be given in the next section.

\begin{lemma}\label{extensions and completions} Let $R$ be an $I$-adically complete noetherian ring, and let $A$ be an abelian scheme over $R$. Then the natural map
\[\operatorname{Ext}^2_{R}(A,\mathbb{G}_m)\to \lim_n \operatorname{Ext}^2_{R/I^n}(A,\mathbb{G}_m)\]
is injective.
\end{lemma}

We are now in position to prove Theorem~\ref{higher vanishing for abelian schemes}.

\begin{proof}[Proof of Theorem~\ref{higher vanishing for abelian schemes}]
Let $T\to S$ be a morphism of schemes, and let $\mathscr{E}\in\operatorname{Ext}^2_T(A,\mathbb{G}_m)$ be an extension class. By Lemma~\ref{Ext is torsion-free}, it suffices to show that $\mathscr{E}$ is fppf-locally torsion. After passing to affine covers, we may assume $T=\operatorname{Spec}R$ and $S=\operatorname{Spec}R_0$. By Corollary~\ref{Ext commutes with filtered colimits} and \cite[Tag \href{https://stacks.math.columbia.edu/tag/00QN}{00QN}]{Stacks}, we further reduce to the case where $R$ is a finitely presented $R_0$-algebra. Finally, using \cite[Rem.\ 27.91]{GWII} and \cite[Tag \href{https://stacks.math.columbia.edu/tag/05N9}{05N9}]{Stacks}, we can assume that $R_0$ is of finite type over $\mathbb{Z}$.

Our goal is to show that for every point $t \in T$, there exists a flat morphism $g\colon T^\prime \to T$ that is locally of finite presentation, whose image contains $t$ and such that the pullback $g^*\mathscr{E}$ is torsion. Let $B\colonequals \mathcal{O}_{T,t}$ be the local ring at $t$, let $\mathfrak{m}$ be its maximal ideal, and consider its $\mathfrak{m}$-adic completion $\widehat{B}$. For each $n\geq 1$, the composition
\[\widehat{B}\to \widehat{B}/\mathfrak{m}^n\to (\widehat{B}/\mathfrak{m}^n)_\text{red}
\simeq \widehat{B}/\mathfrak{m}\]
agrees with the natural quotient map. Lemmas \ref{Ext and reduced schemes} and \ref{extensions and completions} then imply that the maps
\[\operatorname{Ext}^2_{\widehat{B}}(A,\mathbb{G}_m)\to \lim_n \operatorname{Ext}^2_{\widehat{B}/\mathfrak{m}^n}(A,\mathbb{G}_m) \to \prod_{n=1}^\infty \operatorname{Ext}^2_{\widehat{B}/\mathfrak{m}}(A,\mathbb{G}_m)\]
are both injective. According to \cite[\S7]{breen1969extensions}, the extension group $\operatorname{Ext}^2_{\widehat{B}/\mathfrak{m}}(A,\mathbb{G}_m)$ is torsion, and therefore so is $\operatorname{Ext}^2_{\widehat{B}}(A,\mathbb{G}_m)$. (Note that the image of the morphism above lies in the diagonal.)

Given that the local ring $B=\mathcal{O}_{T,t}$ is excellent, Popescu's theorem \cite[Tag \href{https://stacks.math.columbia.edu/tag/07GC}{07GC}]{Stacks} states that the completion $\widehat{B}$ is a filtered colimit of smooth $B$-algebras $B_\lambda$. By Corollary~\ref{Ext commutes with filtered colimits}, the pullback of $\mathscr{E}$ to $\operatorname{Spec}B_\lambda$ is torsion for some $\lambda$, which we now fix.

Now, since $B$ is the filtered colimit of the coordinate rings of affine open neighborhoods of $t$, we may spread out the morphism $\operatorname{Spec}B_\lambda\to \operatorname{Spec}B$ to a smooth map $Y\to V$, where $V$ is an affine neighborhood of $t$. The image of $Y \to V$ contains $t$, because the image of ${\operatorname{Spec}}\widehat{B} \to {\operatorname{Spec}}B\to T$ does. By Corollary~\ref{Ext commutes with filtered colimits}, after possibly shrinking $Y$, the class $\mathscr{E}$ becomes torsion on $Y$, completing the proof.
\end{proof}

\section{Proof of Lemma~\ref{extensions and completions}}\label{section on completions}
Let $R$ be an $I$-adically complete noetherian ring, and denote the quotient $R/I^n$ as $R_n$. This gives rise to a direct system of topoi
\[\textsf{Sh}((\textsf{Sch}/\operatorname{Spec}R_1)_\text{fppf})\xrightarrow{\ f_1\ }\textsf{Sh}((\textsf{Sch}/\operatorname{Spec}R_2)_\text{fppf})\xrightarrow{\ f_2\ }\cdots,\]
whose \emph{lax} colimit exists and will be denoted by $\textsf{X}$ \cite[Thm.\ 2.5]{moerdijk1988classifying}. This colimit coincides with the lax limit of the corresponding diagram of categories with pullback functors:
\[\cdots\xrightarrow{\ f_2^*\ }\textsf{Sh}((\textsf{Sch}/\operatorname{Spec}R_2)_\text{fppf})\xrightarrow{\ f_1^*\ }\textsf{Sh}((\textsf{Sch}/\operatorname{Spec}R_1)_\text{fppf}).\]

Concretely, an object of $\textsf{X}$ is a collection $(F_n,\alpha_n)$, where each $F_n$ is a sheaf on $(\textsf{Sch}/\operatorname{Spec}R_n)_\text{fppf}$ and $\alpha_n\colon f_n^*F_{n+1}\to F_n$ is a morphism on $\textsf{Sh}((\textsf{Sch}/\operatorname{Spec}R_n)_\text{fppf})$. A morphism $(F_n,\alpha_n)\to (G_n,\beta_n)$ is a collection of morphisms $\varphi_n\colon F_n\to G_n$ making the diagram
\[\begin{tikzcd}
	{f_n^* F_{n+1}} && {f_n^* G_{n+1}} \\
	{F_n} && {G_n}
	\arrow["{f_n^*\varphi_{n+1}}", from=1-1, to=1-3]
	\arrow["{\alpha_n}"', from=1-1, to=2-1]
	\arrow["{\beta_n}", from=1-3, to=2-3]
	\arrow["{\varphi_n}", from=2-1, to=2-3]
\end{tikzcd}\]
commute for all $n$. 

The category $\textsf{Ab}(\textsf{X})$ of abelian group objects in $\textsf{X}$ satisfies the universal property of the lax limit
\[\cdots\xrightarrow{\ f_2^*\ }\textsf{Ab}((\textsf{Sch}/\operatorname{Spec}R_2)_\text{fppf})\xrightarrow{\ f_1^*\ }\textsf{Ab}((\textsf{Sch}/\operatorname{Spec}R_1)_\text{fppf})\]
in the 2-category of abelian categories and exact functors, and so can be described in a similar way. An argument, explained to us by Kresch and Mathur and analogous to that of \cite[Prop.\ 1.1]{jannsen1988continuous}, yields a characterization of the injective objects in $\textsf{Ab}(\textsf{X})$.

\begin{lemma}\label{injective}
    Let $(\mathscr{I}_n,\iota_n)$ be an object of $\normalfont\textsf{Ab}(\textsf{X})$, and let $\bar{\iota}_n\colon \mathscr{I}_{n+1}\to f_{n,*}\mathscr{I}_n$ denote the morphism adjoint to $\iota_n$. Then $(\mathscr{I}_n,\iota_n)$ is injective if and only if each $\mathscr{I}_n$ is an injective object of $\normalfont\textsf{Ab}((\textsf{Sch}/\operatorname{Spec}R_n)_\text{fppf})$ and each $\bar{\iota}_n$ is a split epimorphism.
\end{lemma}

\begin{proof}
    Fix a positive integer $m$, and for each $n\leq m$, denote by $f_{n,m}\colon \operatorname{Spec}R_n\to \operatorname{Spec}R_m$ the natural closed immersion. Define a functor
    \[U_m\colon \textsf{Ab}((\textsf{Sch}/\operatorname{Spec}R_m)_\text{fppf})\to \textsf{Ab}(\textsf{X})\]
    by sending an abelian sheaf $\mathscr{G}$ to the object $(\mathscr{G}_n,\gamma_n)$, where
    \[\mathscr{G}_n\colonequals\begin{cases} 
    f_{n,m}^*\mathscr{G} & \text{for } n\leq m \\
   0      & \text{otherwise},
  \end{cases}\]
and the transition map $\gamma_n$ is the natural isomorphism $f_n^*f_{n+1,m}^*\mathscr{G}\simeq f_{n,m}^*\mathscr{G}$ for $n<m$ and zero otherwise.

    Since inverse image functors are exact, so is the functor $U_m$. Moreover, $U_m$ is left adjoint to the functor $V_m\colon \textsf{Ab}(\textsf{X})\to \textsf{Ab}((\textsf{Sch}/\operatorname{Spec}R_m)_\text{fppf})$ that maps $(\mathscr{A}_n,\alpha_n)$ to $\mathscr{A}_m$. By \cite[Tag \href{https://stacks.math.columbia.edu/tag/015Z}{015Z}]{Stacks}, it follows that $V_m$ preserves injectives. In other words, if $(\mathscr{I}_n,\iota_n)$ is an injective object of $\normalfont\textsf{Ab}(\textsf{X})$, then each component $\mathscr{I}_n$ is also injective.

    To show that the adjoint maps $\bar{\iota}_n\colon \mathscr{I}_{n+1}\to f_{n,*}\mathscr{I}_n$ are split epimorphisms, consider the functors
    \[V \colon \textsf{Ab}(\textsf{X}) \to \prod_{n=1}^\infty\textsf{Ab}((\textsf{Sch}/\operatorname{Spec}R_n)_\text{fppf}) \quad\text{and}\quad
    P \colon \prod_{n=1}^\infty\textsf{Ab}((\textsf{Sch}/\operatorname{Spec}R_n)_\text{fppf}) \to \textsf{Ab}(\textsf{X})\]
    defined by $V(\mathscr{A}_n,\alpha_n) \colonequals (\mathscr{A}_n)$ and $P(\mathscr{G}_n) \colonequals \left(\prod_{i=1}^n f_{i,n,*}\mathscr{G}_i,\pi_n\right)$. Here, the transition map $\pi_n$ is defined so that its adjoint
    \[\bar{\pi}_n\colon \prod_{i=1}^{n+1} f_{i,n+1,*}\mathscr{G}_i\to f_{n,*}\left(\prod_{i=1}^n f_{i,n,*}\mathscr{G}_i\right)\simeq \prod_{i=1}^{n} f_{i,n+1,*}\mathscr{G}_i\]
    is the projection onto the first $n$ components.

    The functor $V$ is left adjoint to $P$. Since $V$ is faithful, the unit $\eta\colon (\mathscr{I}_n,\iota_n)\to PV(\mathscr{I}_n,\iota_n)$ is a monomorphism. Using the injectivity of $(\mathscr{I}_n,\iota_n)$, this monomorphism admits a retraction $\nu\colon PV(\mathscr{I}_n,\iota_n)\to (\mathscr{I}_n,\iota_n)$.
\[\begin{tikzcd}
	{\mathscr{I}_{n+1}} && {\displaystyle\prod_{i=1}^{n+1} f_{i,n+1,*}\mathscr{I}_i} && {\mathscr{I}_{n+1}} \\
	{f_{n,*}\mathscr{I}_n} && {f_{n,*}\mleft(\displaystyle\prod_{i=1}^n f_{i,n,*}\mathscr{I}_i\mright)} && {f_{n,*}\mathscr{I}_n}
	\arrow["{\eta_{n+1}}", from=1-1, to=1-3]
	\arrow[bend left=20, equals, from=1-1, to=1-5]
	\arrow["{\bar{\iota}_n}", from=1-1, to=2-1]
	\arrow["{\nu_{n+1}}", from=1-3, to=1-5]
	\arrow["{\bar{\pi}_n}", from=1-3, to=2-3]
	\arrow["{\bar{\iota}_n}", from=1-5, to=2-5]
	\arrow["{f_{n,*}\eta_n}", from=2-1, to=2-3]
	\arrow[bend right=20, equals, from=2-1, to=2-5]
	\arrow["{f_{n,*}\nu_n}", from=2-3, to=2-5]
\end{tikzcd}\]
    Now, let $\rho_n$ be a section of $\bar{\pi}_n$. Then the composition $\nu_{n+1}\circ \rho_n\circ f_{n,*}\eta_n$ is a section of $\bar{\iota}_n$, showing that $\bar{\iota}_n$ is a split epimorphism. Indeed, the commutativity of the diagram above implies that
    \begin{align*}
        \bar{\iota}_n\circ \nu_{n+1}\circ \rho_n\circ f_{n,*}\eta_n &= f_{n,*}\nu_n\circ\bar{\pi}_n\circ \rho_n\circ f_{n,*}\eta_n\\
        &= f_{n,*}\nu_n\circ f_{n,*}\eta_n = \operatorname{id}.
    \end{align*}

    Conversely, let $(\mathscr{I}_n,\iota_n)$ be an object of $\textsf{Ab}(\textsf{X})$ such that each $\mathscr{I}_n$ is an injective object of $\normalfont\textsf{Ab}((\textsf{Sch}/\operatorname{Spec}R_n)_\text{fppf})$ and each $\bar{\iota}_n$ is a split epimorphism. Then $\ker \bar{\iota}_n$ is a direct summand of $\mathscr{I}_{n+1}$, hence injective. Another application of \cite[Tag \href{https://stacks.math.columbia.edu/tag/015Z}{015Z}]{Stacks} shows that $P(\ker\bar{\iota}_n)$ is an injective object of $\textsf{Ab}(\textsf{X})$. But $P(\ker\bar{\iota}_{n-1})$ is isomorphic to $(\mathscr{I}_n,\iota_n)$, concluding the proof.
\end{proof}

Let $(\mathscr{G}_n,\gamma_n)$ and $(\mathscr{A}_n,\alpha_n)$ be objects of $\textsf{Ab}(\textsf{X})$. Suppose that each transition map $\gamma_n\colon f_n^*\mathscr{G}_{n+1}\to \mathscr{G}_n$ is an isomorphism. Given a morphism $\varphi_{n+1}\colon \mathscr{G}_{n+1}\to \mathscr{A}_{n+1}$ in $\textsf{Ab}((\textsf{Sch}/\operatorname{Spec}R_{n+1})_\text{fppf})$, we define 
\[T_n(\varphi_{n+1})\colonequals\left( \mathscr{G}_n \xrightarrow{\ \gamma_n^{-1}\ }f_n^*\mathscr{G}_{n+1}\xrightarrow{\ f_n^*\varphi_{n+1}\ } f_n^*\mathscr{A}_{n+1} \xrightarrow{\ \alpha_n\ }\mathscr{A}_n\right).\]
This yields an inverse system $(\operatorname{Hom}_{R_n}(\mathscr{G}_n,\mathscr{A}_n),T_n)$, and we obtain a natural isomorphism
\[\operatorname{Hom}_\textsf{X}((\mathscr{G}_n,\gamma_n),(\mathscr{A}_n,\alpha_n))\simeq \lim \operatorname{Hom}_{R_n}(\mathscr{G}_n,\mathscr{A}_n).\]
Indeed, an element of $\operatorname{Hom}_\textsf{X}((\mathscr{G}_n,\gamma_n),(\mathscr{A}_n,\alpha_n))$ is a sequence of morphisms $(\varphi_n)$ such that $\varphi_n\circ \gamma_n=\alpha_n\circ f_n^*\varphi_{n+1}$, which is equivalent to
\[\varphi_n=\alpha_n\circ f_n^*\varphi_{n+1}\circ \gamma_n^{-1}=T_n(\varphi_{n+1}).\]

More generally, under the same assumption that each $\gamma_n$ is an isomorphism, we obtain a description of the Ext groups in $\textsf{Ab}(\textsf{X})$.

\begin{proposition}
    Let $(\mathscr{G}_n,\gamma_n)$ and $(\mathscr{A}_n,\alpha_n)$ be objects of $\normalfont\textsf{Ab}(\textsf{X})$, and suppose that each transition map $\gamma_n\colon f_n^*\mathscr{G}_{n+1}\to \mathscr{G}_n$ is an isomorphism. Then, for all $i\geq 0$, there is a short exact sequence
    \[0\to \operatorname{lim}^1\operatorname{Ext}^{i-1}_{R_n}(\mathscr{G}_n,\mathscr{A}_n)\to \operatorname{Ext}^i_{\normalfont\textsf{X}}((\mathscr{G}_n,\gamma_n),(\mathscr{A}_n,\alpha_n))\to \lim \operatorname{Ext}^i_{R_n}(\mathscr{G}_n,\mathscr{A}_n)\to 0,\]
    that is functorial in both arguments.
\end{proposition}

\begin{proof}
The discussion above shows that the functor $\operatorname{Hom}_\textsf{X}((\mathscr{G}_n,\gamma_n),-)\colon \textsf{Ab}(\textsf{X})\to \textsf{Ab}$ is naturally isomorphic to the composition
\[\textsf{Ab}(\textsf{X})\xrightarrow{\ F\ } \textsf{Ab}^\mathbb{N}\xrightarrow{\ \lim \ } \textsf{Ab},\]
where the functor $F$ sends $(\mathscr{A}_n,\alpha_n)$ to the inverse system $(\operatorname{Hom}_{R_n}(\mathscr{G}_n,\mathscr{A}_n),T_n)$. The desired short exact sequences then arise from the Grothendieck spectral sequence associated with this composition. In the remainder of this proof, we verify the existence of this spectral sequence and confirm that it yields the expected result.

The Grothendieck spectral sequence exists under the assumption that $F$ maps injective objects to $\lim$-acyclics. Let $(\mathscr{I}_n,\iota_n)$ be an injective object of $\normalfont\textsf{Ab}(\textsf{X})$, and let $\bar{\iota}_n\colon \mathscr{I}_{n+1}\to f_{n,*}\mathscr{I}_n$ denote the morphism adjoint to $\iota_n$. By Lemma~\ref{injective}, each $\bar{\iota}_n$ admits a section $\sigma_n\colon f_{n,*}\mathscr{I}_n\to \mathscr{I}_{n+1}$. Given any morphism $\varphi_n\colon \mathscr{G}_{n}\to \mathscr{I}_{n}$, define $\varphi_{n+1}$ as the composition
\[\mathscr{G}_{n+1}\xrightarrow{\ \bar{\gamma}_n\ }f_{n,*}\mathscr{G}_n\xrightarrow{\ f_{n,*}\varphi_n\ }f_{n,*}\mathscr{I}_n\xrightarrow{\ \sigma_n\ }\mathscr{I}_{n+1}.\]
It follows that $T_n(\varphi_{n+1})=\varphi_n$, showing that each transition map $T_n$ is surjective. Therefore, the inverse system $(\operatorname{Hom}_{R_n}(\mathscr{G}_n,\mathscr{I}_n),T_n)$ satisfies the Mittag--Leffler condition and is $\lim$-acyclic.

Since $\mathsf{R}^i \lim = 0$ for all $i > 1$ on inverse systems of abelian groups, the Grothendieck spectral sequence degenerates into a family of short exact sequences
\[0\to \operatorname{lim}^1 \textsf{R}^{i-1}F(\mathscr{A}_n,\alpha_n)\to\operatorname{Ext}^i_\textsf{X}((\mathscr{G}_n,\gamma_n),(\mathscr{A}_n,\alpha_n))\to\lim \textsf{R}^{i}F(\mathscr{A}_n,\alpha_n)\to 0.\]
It remains to identify $\textsf{R}^{i}F(\mathscr{A}_n,\alpha_n)$ with the inverse system $(\operatorname{Ext}^i_{R_n}(\mathscr{G}_n,\mathscr{A}_n),T_n^i)$, where the transition maps
\[T_n^i\colon \operatorname{Ext}^i_{R_{n+1}}(\mathscr{G}_{n+1},\mathscr{A}_{n+1})\to \operatorname{Ext}^i_{R_n}(\mathscr{G}_n,\mathscr{A}_n)\]
are induced from the morphisms $T_n$ defined earlier. This follows directly from the description of injective objects in $\textsf{Ab}(\textsf{X})$ given in Lemma~\ref{injective}.
\end{proof}

A natural way to construct objects $(\mathscr{G}_n,\gamma_n)$ with each $\gamma_n$ an isomorphism is as follows. Let $g_n$ denote the closed immersion $\operatorname{Spec}R_n\to \operatorname{Spec}R$. For an abelian sheaf $\mathscr{G}$ on $(\textsf{Sch}/\operatorname{Spec}R)_\text{fppf}$, set $\mathscr{G}_n\colonequals g_n^*\mathscr{G}$. Since $g_n=g_{n+1}\circ f_n$, we obtain natural isomorphisms
\[f_n^*\mathscr{G}_{n+1}= f_n^*g_{n+1}^*\mathscr{G}\simeq g_n^*\mathscr{G}=\mathscr{G}_n,\]
which define the transition maps of an object in $\textsf{Ab}(\textsf{X})$, denoted $\widehat{\mathscr{G}}$.

\begin{corollary}\label{lemma part 1}
   Let $A$ be an abelian scheme over $R$. Then the natural map 
   \[\operatorname{Ext}^2_{\normalfont\textsf{X}}(\widehat{A},\widehat{\mathbb{G}}_m)\to \lim \operatorname{Ext}^2_{R_n}(A,\mathbb{G}_m)\]
   is an isomorphism.
\end{corollary}

\begin{proof}
    From the preceding proposition, we have a short exact sequence
\[0\to \operatorname{lim}^1\operatorname{Ext}^{1}_{R_n}(A,\mathbb{G}_m)\to \operatorname{Ext}^2_{\textsf{X}}(\widehat{A},\widehat{\mathbb{G}}_m)\to \lim \operatorname{Ext}^2_{R_n}(A,\mathbb{G}_m)\to 0.\]
Now, the inverse system $(\operatorname{Ext}^{1}_{R_n}(A,\mathbb{G}_m))$ is Mittag--Leffler since the dual abelian scheme is formally smooth.
\end{proof}

Since the morphisms of topoi $\textsf{Sh}((\textsf{Sch}/\operatorname{Spec}R_n)_\text{fppf})\to\textsf{Sh}((\textsf{Sch}/\operatorname{Spec}R)_\text{fppf})$ factor through $\textsf{X}$, the natural morphism
\[\operatorname{Ext}^i_{R}(\mathscr{G},\mathscr{A})\to\lim \operatorname{Ext}^i_{R_n}(\mathscr{G},\mathscr{A})\]
factors through $\operatorname{Ext}^i_{\textsf{X}}(\widehat{\mathscr{G}},\widehat{\mathscr{A}})$. Having already dealt with $\operatorname{Ext}^i_{\textsf{X}}(\widehat{\mathscr{G}},\widehat{\mathscr{A}})\to \lim \operatorname{Ext}^i_{R_n}(\mathscr{G},\mathscr{A})$, the morphism $\operatorname{Ext}^i_{R}(\mathscr{G},\mathscr{A})\to\operatorname{Ext}^i_{\textsf{X}}(\widehat{\mathscr{G}},\widehat{\mathscr{A}})$ can be studied using the Breen--Deligne resolution.

\begin{lemma}
Let $E$ and $F$ be first-quadrant spectral sequences converging to $E^{p+q}$ and $F^{p+q}$, respectively. Consider a morphism of spectral sequences $f\colon E \to F$ and suppose that, for all $p$, the map $f_1^{p,q}\colon E_1^{p,q} \to F_1^{p,q}$ is an isomorphism for $q \leq 1$ and injective for $q = 2$. Then the induced morphism $E^2 \to F^2$ is injective.
\end{lemma}

\begin{proof}
Recall that $E_2^{p,q}$ is computed as the cohomology of the complex
\[
E_1^{p-1,q} \longrightarrow E_1^{p,q} \longrightarrow E_1^{p+1,q},
\]
and similarly for $F_2^{p,q}$. By the assumption on $f_1^{p,q}$, it follows immediately that $E_2^{p,q} \to F_2^{p,q}$ is an isomorphism for all $p$ when $q \leq 1$, and injective when $(p,q) = (0,2)$.

Let $K_E$ denote the kernel of the natural map $E^2 \to E_2^{0,2}$, and define $K_F$ analogously. We then obtain a commutative diagram with exact rows:
\[
\begin{tikzcd}
0 \arrow{r} & K_E \arrow{r} \arrow{d} & E^2 \arrow{r} \arrow{d} & E_2^{0,2} \arrow[d, hookrightarrow] \\
0 \arrow{r} & K_F \arrow{r} & F^2 \arrow{r} & F_2^{0,2}.
\end{tikzcd}
\]
Since the rightmost vertical arrow is injective, it suffices to show that $K_E \to K_F$ is also injective.

To that end, consider the long exact sequences relating $K_E$ and $K_F$ to adjacent $E_2^{p,q}$ terms. We have a commutative diagram with exact rows:
\[
\begin{tikzcd}
E_2^{0,1} \arrow[d, twoheadrightarrow] \arrow{r} & E_2^{2,0} \arrow[d, hookrightarrow] \arrow{r} & K_E \arrow{d} \arrow{r} & E_2^{1,1} \arrow[d, hookrightarrow] \\
F_2^{0,1} \arrow{r} & F_2^{2,0} \arrow{r} & K_F \arrow{r} & F_2^{1,1},
\end{tikzcd}
\]
in which the indicated vertical arrows are either surjective or injective by our assumptions on $f_1^{p,q}$. A diagram chase now shows that the map $K_E \to K_F$ is injective, completing the proof.
\end{proof}

\begin{corollary}\label{lemma part 2}
   Let $A$ be an abelian scheme over $R$. Then the natural map 
   \[\operatorname{Ext}^2_{R}(A,\mathbb{G}_m)\to \operatorname{Ext}^2_{\normalfont\textsf{X}}(\widehat{A},\widehat{\mathbb{G}}_m)\]
   is injective.
\end{corollary}

\begin{proof}
The preceding lemma, when applied to the Breen--Deligne spectral sequence, shows that the morphism $\operatorname{Ext}^2_{R}(A,\mathbb{G}_m)\to\operatorname{Ext}^2_{\textsf{X}}(\widehat{A},\widehat{\mathbb{G}}_m)$ is injective provided that, for all $n\geq 0$, the comparison map
\[\textrm{H}^i(A^n,\mathbb{G}_m)\to \textrm{H}^i(\widehat{A}^n,\widehat{\mathbb{G}}_m)\]
is an isomorphism for $i\leq 1$ and injective for $i=2$. Here, the object $\textrm{H}^i(\widehat{A}^n,\widehat{\mathbb{G}}_m)$ refers to the cohomology group computed in the topos $\textsf{X}$, as defined in the beginning of Section~\ref{Sect 2}. Since the abelian scheme $A$ was arbitrary, we may assume that $n=1$.

Unwinding the definitions, we observe that for $i=0$, the comparison map is simply the identity on $R^\times$. For $i = 1$ and $i= 2$, it coincides\footnote{Although the cohomology groups in \cite{kresch2023formal} are computed on the small étale sites, the interpretation of $\textrm{H}^1$ and $\textrm{H}^2$ in terms of torsors and gerbes makes it clear that these groups coincide with ours.} with the morphisms considered in \cite[Props.\ 3.2 and 5.3]{kresch2023formal}, thereby completing the proof.
\end{proof}

The Lemma~\ref{extensions and completions} has now been proven as a concatenation of Corollaries~\ref{lemma part 1} and \ref{lemma part 2}.

\section{Extensions of the additive group}
As noted in the introduction, the sheaf $\underline{\operatorname{Ext}}^1(\mathbb{G}_a,\mathbb{G}_m)$ on the site $(\textsf{Sch}/\operatorname{Spec}\mathbb{Q})_\text{fppf}$ has often been incorrectly claimed to vanish in the literature. In this section, we compute it as explicitly as possible.

\begin{theorem}\label{Ext additive}
Let $T$ be a quasi-compact and quasi-separated $\mathbb{Q}$-scheme. Then the natural maps
\[\normalfont\underline{\operatorname{Ext}}^1(\mathbb{G}_a,\mathbb{G}_m)(T)\to \underline{\operatorname{Ext}}^1(\mathbb{G}_a,\mathbb{G}_m)(T_\text{red})\leftarrow \operatorname{Ext}^1_{T_\text{red}}(\mathbb{G}_a,\mathbb{G}_m)\to \textrm{H}^1_m(\mathbb{G}_{a,T_\text{red}},\mathbb{G}_{m,T_\text{red}})\]
are all isomorphisms. These groups vanish if $\normalfont T_\text{red}$ is seminormal. Conversely, if $\normalfont T_\text{red}$ is affine and not seminormal, they are nonzero.
\end{theorem}

\begin{remark}
    Let $h\colon \operatorname{Spec}\mathbb{Q} \to \operatorname{Spec}\mathbb{Z}$ denote the natural morphism. For a scheme $S$, let $\underline{\operatorname{Ext}}^1(\mathbb{G}_{a,S}, \mathbb{G}_{m,S})$ denote the extension sheaf computed in the category $\textsf{Ab}((\textsf{Sch}/S)_\text{fppf})$. In \cite{gabberstalk}, Gabber claimed that the unit map
\[\underline{\operatorname{Ext}}^1(\mathbb{G}_{a,\mathbb{Z}}, \mathbb{G}_{m,\mathbb{Z}}) \to h_* h^* \underline{\operatorname{Ext}}^1(\mathbb{G}_{a,\mathbb{Z}}, \mathbb{G}_{m,\mathbb{Z}}) \simeq h_* \underline{\operatorname{Ext}}^1(\mathbb{G}_{a,\mathbb{Q}}, \mathbb{G}_{m,\mathbb{Q}})\]
is an isomorphism. In other words, for any scheme $T$, we have a natural isomorphism
\[\underline{\operatorname{Ext}}^1(\mathbb{G}_{a,\mathbb{Z}}, \mathbb{G}_{m,\mathbb{Z}})(T) \simeq \underline{\operatorname{Ext}}^1(\mathbb{G}_{a,\mathbb{Q}}, \mathbb{G}_{m,\mathbb{Q}})(T \times \operatorname{Spec}\mathbb{Q}),\]
where the right-hand side can be computed using Theorem~\ref{Ext additive}.
\end{remark}

\begin{remark}
    Since multiple non-equivalent definitions of seminormal rings and schemes appear in the literature, we clarify that this paper adopts Swan’s definition \cite{swan1980seminormality}. A ring $R$ is said to be \emph{seminormal} if, for all $x, y \in R$ satisfying $x^3 = y^2$, there exists $r \in R$ such that $x = r^2$ and $y = r^3$. A scheme $S$ is then called \emph{seminormal} if every point of $S$ admits an affine open neighborhood isomorphic to the spectrum of a seminormal ring. For further details, see \cite[Tag~\href{https://stacks.math.columbia.edu/tag/0EUK}{0EUK}]{Stacks}.
\end{remark}

The proof of Theorem~\ref{Ext additive} will require some auxiliary lemmas, which we now establish.

\begin{lemma}\label{lemma Ga red}
   Let $R$ be a $\mathbb{Q}$-algebra. Then the restriction map
    \[\normalfont\operatorname{Ext}^1_R(\mathbb{G}_a,\mathbb{G}_m)\to \operatorname{Ext}^1_{R_\text{red}}(\mathbb{G}_a,\mathbb{G}_m)\]
    is an isomorphism.
\end{lemma}

\begin{proof}
 According to Corollaries~\ref{H2} and \ref{H3}, the cohomology groups $\textrm{H}^2_s(\mathbb{G}_{a,R}, \mathbb{G}_{m,R})$ and $\textrm{H}^3_s(\mathbb{G}_{a,R}, \mathbb{G}_{m,R})$ vanish---and the same holds for the corresponding groups over $R_{\text{red}}$. Consequently, Proposition~\ref{short exact sequence computing extensions} yields the following commutative diagram:
\[\begin{tikzcd}
	{\operatorname{Ext}^1_R(\mathbb{G}_a,\mathbb{G}_m)} & {\textrm{H}^1_m(\mathbb{G}_{a,R},\mathbb{G}_{m,R})} \\
	{\operatorname{Ext}^1_{R_\text{red}}(\mathbb{G}_a,\mathbb{G}_m)} & {\textrm{H}^1_m(\mathbb{G}_{a,R_\text{red}},\mathbb{G}_{m,R_\text{red}}).}
	\arrow["\sim", from=1-1, to=1-2]
	\arrow[from=1-1, to=2-1]
	\arrow[from=1-2, to=2-2]
	\arrow["\sim", from=2-1, to=2-2]
\end{tikzcd}\]
The conclusion then follows from the classical fact that, for an affine scheme $S$, the restriction map $\operatorname{Pic}(S) \to \operatorname{Pic}(S_{\text{red}})$ is an isomorphism \cite[Lem.\ 2.2.9]{rosengarten2023tate}.
\end{proof}

\begin{lemma}\label{seminormal vanishing}
    Let $T$ be a scheme. If $T$ is seminormal, the group $\normalfont\textrm{H}^1_m(\mathbb{G}_{a,T},\mathbb{G}_{m,T})$ vanishes. Conversely, if $T$ is affine and $\normalfont\textrm{H}^1_m(\mathbb{G}_{a,T},\mathbb{G}_{m,T})$ vanishes, then $\normalfont T_\text{red}$ is seminormal.
\end{lemma}

\begin{proof}
   Let $p\colon \mathbb{G}_{a,T} \to T$ be the structure morphism. If $T$ is seminormal, the pullback map
   \[p^*\colon \operatorname{Pic}(T)\to \operatorname{Pic}(\mathbb{G}_{a,T})\]
   is an isomorphism \cite[Lem.\ 4.3]{sadhu2021equivariant}. It follows that $\textrm{H}^1_m(\mathbb{G}_{a,T}, \mathbb{G}_{m,T})$ identifies with the subgroup of $\operatorname{Pic}(T)$ consisting of those line bundles $x$ for which
   \[m^*p^* x=\operatorname{pr}_1^*p^*x+\operatorname{pr}_2^*p^*x,\]
   where $m\colon \mathbb{G}_{a,T} \times_T \mathbb{G}_{a,T} \to \mathbb{G}_{a,T}$ is the group law and $\operatorname{pr}_i$ are the natural projections. However, the morphisms $p \circ m$, $p \circ \operatorname{pr}_1$, and $p \circ \operatorname{pr}_2$ all agree with the structure morphism of $\mathbb{G}_{a,T}^2$, which admits a section $T \to \mathbb{G}_{a,T}^2$. Therefore, the identity above holds if and only if $x = 0$, showing that $\textrm{H}^1_m(\mathbb{G}_{a,T}, \mathbb{G}_{m,T})$ vanishes.

   Now suppose that $T = \operatorname{Spec}R$ is such that $R_\text{red}$ is not seminormal. We adapt a construction of Schanuel \cite{bass1962torsion} to exhibit a nonzero element of $\textrm{H}^1_m(\mathbb{G}_{a,R}, \mathbb{G}_{m,R})$. As in the proof of Lemma~\ref{lemma Ga red}, there is an isomorphism $\textrm{H}^1_m(\mathbb{G}_{a,R}, \mathbb{G}_{m,R}) \simeq \textrm{H}^1_m(\mathbb{G}_{a,R_\text{red}}, \mathbb{G}_{m,R_\text{red}})$, so we may assume without loss of generality that $R$ is reduced. Let $R^{\text{sn}} \supset R$ denote the seminormalization of $R$.
   
   By \cite[Lem.\ 2.6]{swan1980seminormality}, there exists an element $s \in R^{\text{sn}} \setminus R$ such that $s^2, s^3 \in R$. Consider the $R[t]$-submodules $M=(s^2t^2,1+st)$ and $N=(s^2t^2,1-st)$ of $R^\text{sn}[t]$. Then the tensor product $M\otimes_{R[t]}N$ is generated by the elements
   \[\left\{s^4t^4,s^2t^2+s^3t^3,s^2t^2-s^3t^3,1-s^2t^2\right\}\subset R[t],\]
   and thus naturally lies inside $R[t]$. Moreover, we have the identity
   \[1=s^4t^4+(1-s^2t^2)(1+s^2t^2)\in M\otimes_{R[t]}N,\]
   proving that $M$ is an invertible $R[t]$-module, with inverse $N$. In other words, $M$ defines a nonzero element of $\operatorname{Pic}(\mathbb{G}_{a,R})$.

   We now prove that $M$ belongs to the subgroup $\textrm{H}^1_m(\mathbb{G}_{a,R}, \mathbb{G}_{m,R})$ of $\operatorname{Pic}(\mathbb{G}_{a,R})$. Consider the tensor product $m^*M\otimes_{R[x,y]}\operatorname{pr}_1^* N\otimes_{R[x,y]}\operatorname{pr}_2^* N$, which is generated by
   \begin{gather*}
       \left\{
   s^6(x+y)^2x^2y^2, x^2(s^4-s^5y)(x+y)^2, y^2(s^4-s^5x)(x+y)^2,\right.\\
   (s^2-s^3x-s^3y+s^4xy)(x+y)^2, x^2y^2(s^4+s^5x+s^5y), x^2(s^2+s^3x-s^4y^2-s^4xy), \\ \left. y^2(s^2+s^3y-s^4x^2-s^4xy), 1-s^2x^2-s^2y^2-s^3x^2y-s^2xy-s^3xy^2 \right\}\subset R[x,y].
   \end{gather*}
   It follows that it lies inside $R[x,y]$. As one can verify in their preferred programming language, this submodule also contains $1$, proving that $M$ is multiplicative.
\end{proof}

\begin{remark}
    For the reader’s convenience, we include a Macaulay2 implementation of an algorithm that verifies the submodule $m^*M\otimes_{R[x,y]}\operatorname{pr}_1^* N\otimes_{R[x,y]}\operatorname{pr}_2^* N$ of $R[x,y]$ in the preceding proof contains the unit $1$.
       	\begin{lstlisting}[language=Macaulay2,extendedchars=true]
 -- Define the ring ZZ[x,y,s^2,s^3] and the ideal 'prod' with the given generators
 Rxy = ZZ[x, y, u, v] / (v^2 - u^3)
 prod = ideal(
        v^2*(x + y)^2*x^2*y^2,
        x^2*(u^2 - u*v*y)*(x + y)^2,
        y^2*(u^2 - u*v*x)*(x + y)^2,
        (u - v*x - v*y + u^2*x*y)*(x + y)^2,
        x^2*y^2*(u^2 + u*v*x + u*v*y),
        x^2*(u + v*x - u^2*y^2 - u^2*x*y),
        y^2*(u + v*y - u^2*x^2 - u^2*x*y),
        1 - u*x^2 - u*y^2 - v*x^2*y - u*x*y - v*x*y^2
 )

 -- Compute a Groebner basis and the change of basis matrix, which provides the coefficients for writing 1 as a linear combination of the generators of 'prod'
 grob = gb(prod, ChangeMatrix => true)
 mat = getChangeMatrix(grob)

 -- Print the result
 gens grob
\end{lstlisting}
\end{remark}

The lemma above shows that for an \emph{affine} scheme $T$, the group $\textrm{H}^1_m(\mathbb{G}_{a,T},\mathbb{G}_{m,T})$ vanishes if and only if $T_\text{red}$ is seminormal. We now present an example, due to Weibel, which demonstrates that this characterization fails for non-affine schemes.

\begin{remark}\label{weibels example}
Let $k$ be a field, and let $\mathscr{L} \colonequals \mathcal{O}(-1)$ denote the tautological line bundle on the projective line $\mathbb{P}^1_k$. Consider its symmetric algebra:
\[\operatorname{Sym}\mathscr{L}=\bigoplus_{n=0}^\infty \operatorname{Sym}^n\mathscr{L}.\]
By omitting the degree-one component, we obtain a graded subalgebra $\mathscr{A} \subset \operatorname{Sym} \mathscr{L}$, defined as
\[\mathscr{A}\colonequals \operatorname{Sym}^0\mathscr{L} \oplus \operatorname{Sym}^2\mathscr{L}\oplus \operatorname{Sym}^3\mathscr{L}\oplus \dots.\]

For an affine open subset $W = \operatorname{Spec} R \subset \mathbb{P}^1_k$ with $\operatorname{Pic}(W)=0$ (for instance, the complement of a $k$-rational point), the algebra $\mathscr{A}(W)$ is isomorphic to the cusp ring $R[t^2, t^3]$, where $t$ is a generator of the free $R$-module $\mathscr{L}(W)$. This ring is a prototypical example of a reduced but non-seminormal ring. Consequently, the relative spectrum $T$ of $\mathscr{A}$ is a reduced scheme that is not seminormal.

The natural inclusion $R[t^2,t^3]\hookrightarrow R[t]$, together with the ideal $(t^2,t^3)\subset R[t^2,t^3]$, gives rise to a Milnor square\footnote{See \cite[\S I.2.6]{KBook} for the definition.}
\[\begin{tikzcd}[ampersand replacement=\&]
	{R[t^2,t^3]} \& {R[t]} \\
	R \& {R[t]/(t^2),}
	\arrow[hook, from=1-1, to=1-2]
	\arrow[from=1-1, to=2-1]
	\arrow[two heads, from=1-2, to=2-2]
	\arrow[hook, from=2-1, to=2-2]
\end{tikzcd}\]
where the map $R[t^2,t^3]\to R$ sends a polynomial $p(t)$ to its constant term $p(0)$. The associated Units-Pic exact sequence \cite[Thm.\ I.3.10]{KBook} then reads as
\[\begin{tikzcd}
	0 & {R[t^2,t^3]^\times} & {R[t]^\times\times R^\times} & {\left(R[t]/(t^2)\right)^\times}\ar[draw=none]{dll}[name=X, anchor=center]{}\ar[rounded corners,
            to path={ -- ([xshift=2ex]\tikztostart.east)
                      |- (X.center) \tikztonodes
                      -| ([xshift=-2ex]\tikztotarget.west)
                      -- (\tikztotarget)}]{dll}[at end]{} \\
	& {\operatorname{Pic}(R[t^2,t^3])} & {\operatorname{Pic}(R[t])\times\operatorname{Pic}(R)} & {\operatorname{Pic}\mleft(R[t]/(t^2)\mright).}
	\arrow[from=1-1, to=1-2]
	\arrow[from=1-2, to=1-3]
	\arrow[from=1-3, to=1-4]
	\arrow[from=2-2, to=2-3]
    \arrow[from=2-3, to=2-4]
\end{tikzcd}\]

Since $R$ is reduced, both $R[t^2,t^3]^\times$ and $R[t]^\times$ equal $R^\times$. Furthermore, the Picard groups $\operatorname{Pic}(R[t])$ and $\operatorname{Pic}(R)$ vanish. The sequence then simplifies to:
\[\begin{tikzcd}[ampersand replacement=\&,row sep=tiny]
	0 \& {R^\times} \& {R^\times\times R^\times} \& {\left(R[t]/(t^2)\right)^\times} \& {\operatorname{Pic}(R[t^2,t^3])} \& 0 \\
	\& r \& {(r,r)} \\
	\&\& {(r,s)} \& {rs^{-1},}
	\arrow[from=1-1, to=1-2]
	\arrow[from=1-2, to=1-3]
	\arrow[from=1-3, to=1-4]
	\arrow[from=1-4, to=1-5]
	\arrow[from=1-5, to=1-6]
	\arrow[maps to, from=2-2, to=2-3]
	\arrow[maps to, from=3-3, to=3-4]
\end{tikzcd}\]
from which we deduce an isomorphism $\operatorname{Pic}(R[t^2,t^3])\simeq (R,+)$.

For a scheme $S$, let $e\colon S\to \mathbb{A}^1_S$ be the zero-section. We denote by $\operatorname{NPic}(S)$ the kernel of the pullback map
\[e^*\colon \operatorname{Pic}(\mathbb{A}^1_S)\to \operatorname{Pic}(S).\]
Note that $\textrm{H}^1_m(\mathbb{G}_{a,S},\mathbb{G}_{m,S})$ is a subgroup of $\operatorname{NPic}(S)$. The same argument as above applies with $R$ replaced by the polynomial ring $R[z]$, yielding the identification
\[\operatorname{NPic}(R[t^2,t^3])\simeq zR[z]\simeq R[z].\]

Let $U=\operatorname{Spec}k[t]$ and $V=\operatorname{Spec}k[t^{-1}]$ be the standard open cover of $\mathbb{P}^1_k$. Since $\operatorname{NPic}$ is a Zariski sheaf on reduced schemes \cite[Thm.\ 4.7]{weibel1991pic}, the group $\operatorname{NPic}(T)$ is the equalizer of the diagram
\[\operatorname{NPic}(X_U)\times \operatorname{NPic}(X_{V}) \simeq k[t,z]\times k[t^{-1},z]\rightrightarrows k[t,t^{-1},z]\simeq \operatorname{NPic}(X_{U\cap V}),\]
where the maps are given by
\[\begin{tikzcd}[ampersand replacement=\&,row sep=small]
	\& {tp(t,z)} \\
	{(p(t,z),q(t^{-1},z))} \\
	\& {q(t^{-1},z).}
	\arrow[maps to, from=2-1, to=1-2]
	\arrow[maps to, from=2-1, to=3-2]
\end{tikzcd}\]
As in the usual computation of $\mathscr{L}(\mathbb{P}^1_k)=0$, we obtain that $\operatorname{NPic}(T)$ vanishes and so does $\textrm{H}^1_m(\mathbb{G}_{a,T},\mathbb{G}_{m,T})$.
\end{remark}

\begin{proof}[Proof of Theorem~\ref{Ext additive}]
We first show that the natural map from $\underline{\operatorname{Ext}}^1(\mathbb{G}_a,\mathbb{G}_m)(T)$ to $\underline{\operatorname{Ext}}^1(\mathbb{G}_a,\mathbb{G}_m)(T_\text{red})$ is an isomorphism. Since this a local question on $T$, we may assume $T=\operatorname{Spec}R$ is affine. By Lemma~\ref{lemma Ga red}, it suffices to show that the vertical maps in the diagram
\[\begin{tikzcd}
	{\operatorname{Ext}^1_R(\mathbb{G}_a,\mathbb{G}_m)} & {\operatorname{Ext}^1_{R_\text{red}}(\mathbb{G}_a,\mathbb{G}_m)} \\
	{\underline{\operatorname{Ext}}^1(\mathbb{G}_a,\mathbb{G}_m)(T)} & {\underline{\operatorname{Ext}}^1(\mathbb{G}_a,\mathbb{G}_m)(T_\text{red})}
	\arrow[from=1-1, to=1-2]
	\arrow[from=1-1, to=2-1]
	\arrow[from=1-2, to=2-2]
	\arrow[from=2-1, to=2-2]
\end{tikzcd}\]
are isomorphisms. 

The sheaf $\underline{\operatorname{Hom}}(\mathbb{G}_a,\mathbb{G}_m)$ is represented by the formal completion $\widehat{\mathbb{G}}_a$ of $\mathbb{G}_a$ along the identity section. By Proposition~\ref{sheafification map}, the vertical maps in the diagram above are isomorphisms provided that $\textrm{H}^n(R,\widehat{\mathbb{G}}_a)$ vanishes for $n=1,2$. This holds due to a computation of de Jong, explained in \cite[Rem.\ 2.2.18]{bhatt2022prismatic}.

We now assume that $T$ is reduced and prove that $\operatorname{Ext}^1_{T}(\mathbb{G}_a,\mathbb{G}_m)\to \underline{\operatorname{Ext}}^1(\mathbb{G}_a,\mathbb{G}_m)(T)$ and  $\operatorname{Ext}^1_{T}(\mathbb{G}_a,\mathbb{G}_m)\to \textrm{H}^1_m(\mathbb{G}_{a,T},\mathbb{G}_{m,T})$ are isomorphisms. As above, the first map is an isomorphism as soon as $\textrm{H}^n(T,\widehat{\mathbb{G}}_a)$ vanishes for $n=1,2$. To show this, we choose a Zariski hypercover $K=(I,\{U_i\})$ of $T$ such that each $U_i$ is an affine open of $T$, whose existence is guaranteed by \cite[Tag \href{https://stacks.math.columbia.edu/tag/01H7}{01H7}]{Stacks}. Then, \cite[Tag \href{https://stacks.math.columbia.edu/tag/01GY}{01GY}]{Stacks} gives a spectral sequence
\[E_2^{p,q}\implies \textrm{H}^{p+q}(T,\widehat{\mathbb{G}}_a),\]
where $E_2^{p,q}$ is the $p$-th cohomology group of the complex
\[\prod_{i\in I_0}\textrm{H}^{q}(U_i,\widehat{\mathbb{G}}_a)\to \prod_{i\in I_1}\textrm{H}^{q}(U_i,\widehat{\mathbb{G}}_a)\to \prod_{i\in I_2}\textrm{H}^{q}(U_i,\widehat{\mathbb{G}}_a)\to \dots.\]

Since each $U_i$ is affine, de Jong's computation implies that $E_2^{p,q}=0$ for all $p$ and all $q>0$. As a result, the spectral sequence degenerates and we have that $\textrm{H}^n(T,\widehat{\mathbb{G}}_a)$ is the $n$-th cohomology group of the complex
\[\prod_{i\in I_0}\Gamma(U_i,\widehat{\mathbb{G}}_a)\to \prod_{i\in I_1}\Gamma(U_i,\widehat{\mathbb{G}}_a)\to \prod_{i\in I_2}\Gamma(U_i,\widehat{\mathbb{G}}_a)\to \dots.\]
By \cite[Lem.\ 2.2.5]{rosengarten2023tate}, the fact that each $U_i$ is the spectrum of a reduced ring implies that $\Gamma(U_i,\widehat{\mathbb{G}}_a)=0$, and hence $\textrm{H}^n(T,\widehat{\mathbb{G}}_a)$ vanishes for all $n$.

Next, Proposition~\ref{short exact sequence computing extensions} shows that the map $\operatorname{Ext}^1_{T}(\mathbb{G}_a,\mathbb{G}_m)\to \textrm{H}^1_m(\mathbb{G}_{a,T},\mathbb{G}_{m,T})$ is an isomorphism if $\textrm{H}^n_s(\mathbb{G}_{a,T},\mathbb{G}_{m,T})$ vanishes for $n=2,3$. This vanishing follows from the fact that, since $T$ is assumed reduced, any morphism of $T$-schemes $\mathbb{G}_{a,T}^n\to \mathbb{G}_{m,T}$ must be constant. To complete the proof, it remains to show that $\textrm{H}^1_m(\mathbb{G}_{a,T},\mathbb{G}_{m,T})$ vanishes if $T$ is seminormal and that it does not vanish when $T$ is affine and not seminormal. Both facts were established in Lemma~\ref{seminormal vanishing}.
\end{proof}

\section{Group cohomology computations}
In this section, we gather several group cohomology computations used for the application of Proposition~\ref{short exact sequence computing extensions} in the previous section.

\begin{proposition}[Lazard]\label{Lazard's theorem}
Let $R$ be a ring, and let $S$ denote the set of prime numbers that are not invertible in $R$. Then, we have an isomorphism of $R$-modules
\[\normalfont\textrm{H}^2_s(\mathbb{G}_{a,R}, \mathbb{G}_{a,R})\simeq \bigoplus_{\substack{p\in S\\ n\geq 1}}\left[Q_{p^n}(x,y)\right]\cdot R/pR,\]
for $Q_{p^n}(x,y)\colonequals ((x+y)^{p^n} - x^{p^n} - y^{p^n})/p$. Moreover, $\normalfont\textrm{H}^2_s(\mathbb{G}_{a,R}, \mathbb{G}_{a,R})$ coincides with the full second Hochschild cohomology group $\normalfont\textrm{H}^2_0(\mathbb{G}_{a,R}, \mathbb{G}_{a,R})$ if and only if $R$ is torsion-free as an abelian group.
\end{proposition}

\begin{proof}
According to Definition~\ref{def invariants}, $\textrm{H}^2_s(\mathbb{G}_{a,R}, \mathbb{G}_{a,R})$ is the middle cohomology of the complex of $R$-modules
\[R[x]\xrightarrow{\ \delta\ } R[x,y]\xrightarrow{\ \gamma\ } R[x,y,z]\oplus R[x,y],\]
whose differentials are given by
\begin{align*}
    \delta(p(x)) &= p(x+y)-p(x)-p(y)\\
    \gamma(q(x,y)) &= (q(x+y,z)+q(x,y)-q(x,y+z)-q(y,z),q(x,y)-q(y,x)).
\end{align*}

In \cite[Lem.\ 3]{Lazard}, Lazard proves that the only homogeneous polynomials of degree $d\geq 1$ in the kernel of $\gamma$ are scalar multiples of
\[Q_d(x,y)\colonequals \begin{cases} 
   \frac{1}{p}\left((x+y)^d - x^d - y^d\right) & \text{if } d \text{ is a power of a prime number }p \\
   (x+y)^d - x^d - y^d       & \text{otherwise}.
  \end{cases}\]
When $d$ is a power of $p$, the polynomial $Q_d(x,y)$ is in the image of $\delta$ if and only if $p\in R^\times$. When $d$ is not a prime power, $Q_d(x,y)$ always lies in the image of $\delta$. Constant polynomials lie in both the kernel of $\gamma$ and in the image of $\delta$, hence do not contribute to cohomology.

As remarked by Lazard in \cite[Rem.\ 3.14]{Lazard}, if $R$ is torsion-free as an abelian group, then every polynomial in the kernel of the map
\begin{align*}
    \widetilde{\gamma}\colon R[x,y] &\to R[x,y,z]\\
    q(x,y) &\mapsto q(x+y,z)+q(x,y)-q(x,y+z)-q(y,z)
\end{align*}
is automatically symmetric. Hence, we have that $\textrm{H}^2_s(\mathbb{G}_{a,R}, \mathbb{G}_{a,R})=\textrm{H}^2_0(\mathbb{G}_{a,R}, \mathbb{G}_{a,R})$. 

For the converse, suppose that $R$ is not torsion-free as an abelian group and let $n$ be the smallest positive integer such that there exists $r\in R\setminus\{0\}$ with $nr=0$. Pick a prime $p$ dividing $n$, and set $k=n/p$ so that $s\colonequals kr\neq 0$ and $ps=0$. Consider the polynomial $sx^py\in R[x,y]$. Then
\[\widetilde{\gamma}(sx^py)=s(x+y)^pz-sx^pz-sy^pz=\sum_{i=1}^{p-1}\binom{p}{i}sx^iy^{p-i}z\]
vanishes because each $\binom{p}{i}$, for $0<i<p$, is divisible by $p$. As a result, $sx^py$ defines a class in $\textrm{H}^2_0(\mathbb{G}_{a,R}, \mathbb{G}_{a,R})$ which does not lie in $\textrm{H}^2_s(\mathbb{G}_{a,R}, \mathbb{G}_{a,R})$.
\end{proof}

\begin{corollary}\label{H2}
Let $R$ be a ring. If either $R$ is reduced or a $\mathbb{Q}$-algebra, then the second Hochschild cohomology group $\normalfont\textrm{H}^2_0(\mathbb{G}_{a,R}, \mathbb{G}_{m,R})$ vanishes. In particular, so does $\normalfont\textrm{H}^2_s(\mathbb{G}_{a,R}, \mathbb{G}_{m,R})$.
\end{corollary}

\begin{proof}
    By definition, the abelian group $\textrm{H}^2_0(\mathbb{G}_{a,R}, \mathbb{G}_{m,R})$ is computed as the cohomology of the complex
    \[R[x]^\times\xrightarrow{\ \delta^\prime\ } R[x,y]^\times\xrightarrow{\ \gamma^\prime\ } R[x,y,z]^\times,\]
with differentials given by
\begin{align*}
    \delta^\prime(p(x)) &= p(x+y)/p(x)p(y)\\
    \gamma^\prime(q(x,y)) &= q(x+y,z)q(x,y)/q(x,y+z)q(y,z).
\end{align*}
This complex is exact when $R$ is reduced, since the units in the corresponding polynomial rings are then necessarily constant. We therefore assume from now on that $R$ is a $\mathbb{Q}$-algebra.

Let $q(x,y)\in R[x,y]^\times$ be an element in the kernel of $\gamma^\prime$. By Remark~\ref{normalized BD}, we may assume that $q(0,0)=1$ and write $q(x,y)=1+n(x,y)$, where $n(x,y)$ is a polynomial with nilpotent coefficients. Since $R$ is a $\mathbb{Q}$-algebra, we have that
\[q(x,y)=\operatorname{exp}(g(x,y)), \quad \text{for}\quad g(x,y)=\log(q(x,y))=\sum_{k=1}^\infty \frac{(-1)^{k+1}}{k}n(x,y)^k.\]
The polynomial $g(x,y)$ lies in the kernel of the map $\widetilde{\gamma}$ from the proof of the previous proposition. By its main result, there exists a polynomial $p(x)\in R[x]$, with nilpotent coefficients, such that $\delta(p(x))=g(x,y)$, where $\delta$ is the additive coboundary defined there. It then follows that $f(x)\colonequals \operatorname{exp}(p(x))$ satisfies $\delta^\prime(f(x))=q(x,y)$. Hence, $q(x,y)$ lies in the image of $\delta^\prime$, completing the proof.
\end{proof}

\begin{remark}
Let $R$ be a reduced ring, and let $R^\prime$ denote its ring of dual numbers $R[\varepsilon]/(\varepsilon^2)$. Then any element $p(x)\in R^\prime[x]^\times$ satisfying $p(0)=1$ is of the form $1+\varepsilon f(x)$, for a unique polynomial $f(x)\in R[x]$ with $f(0)=0$. The same holds for polynomial rings in more variables. It follows that
\[\textrm{H}^2_s(\mathbb{G}_{a,R^\prime}, \mathbb{G}_{m,R^\prime})\simeq \textrm{H}^2_s(\mathbb{G}_{a,R}, \mathbb{G}_{a,R}),\]
proving that $\textrm{H}^2_s(\mathbb{G}_{a,S}, \mathbb{G}_{m,S})$ may be nonzero if $S$ is a ring that is not reduced nor a $\mathbb{Q}$-algebra.
\end{remark}

Let $R$ be a ring. Just as the group $\textrm{H}^2_s(\mathbb{G}_{a,R}, \mathbb{G}_{a,R})$ naturally embeds into $\textrm{H}^2_0(\mathbb{G}_{a,R}, \mathbb{G}_{a,R})$, there is a morphism
\[\eta\colon \textrm{H}^3_s(\mathbb{G}_{a,R}, \mathbb{G}_{a,R})\to \textrm{H}^3_0(\mathbb{G}_{a,R}, \mathbb{G}_{a,R}),\]
induced by the morphism of complexes
\[\begin{tikzcd}[ampersand replacement=\&]
	{R[x,y]} \& {R[x,y,z]\oplus R[x,y]} \& {R[x,y,z,w]\oplus R[x,y,z]^{\oplus 2}\oplus R[x,y]\oplus R[x]} \\
	{R[x,y]} \& {R[x,y,z]} \& {R[x,y,z,w],}
	\arrow["\gamma", from=1-1, to=1-2]
	\arrow[equals, from=1-1, to=2-1]
	\arrow["\beta", from=1-2, to=1-3]
	\arrow[from=1-2, to=2-2]
	\arrow[from=1-3, to=2-3]
	\arrow["{\widetilde{\gamma}}", from=2-1, to=2-2]
	\arrow["{\widetilde{\beta}}", from=2-2, to=2-3]
\end{tikzcd}\]
where the vertical arrows are the projections onto the first summand.

\begin{proposition}\label{H3 Ga}
Let $R$ be a ring that is torsion-free as an abelian group. Then the natural map
\[\eta\colon \normalfont\textrm{H}^3_s(\mathbb{G}_{a,R}, \mathbb{G}_{a,R})\to \textrm{H}^3_0(\mathbb{G}_{a,R}, \mathbb{G}_{a,R})\]
is injective. Moreover, if $R$ is a $\mathbb{Q}$-algebra, then both cohomology groups vanish.
\end{proposition}

\begin{proof}
Using the notation in the commutative diagram above, let $(p,q)\in \ker \beta$.  That is, the polynomials $p$ and $q$ satisfy the following relations:
    \begin{align} 
   \label{eq 1.1} p(x+y,z,w)+p(x,y,z+w) &= p(y,z,w)+p(x,y+z,w)+p(x,y,z) \\
   \label{eq 1.2} p(x,y,z)+p(z,x,y)+q(y,z)+q(x,z) &= p(x,z,y)+q(x+y,z) \\
   \label{eq 1.3} p(x,y,z)+p(y,z,x)+q(x,y+z) &= p(y,x,z)+q(x,y)+q(x,z) \\
   \label{eq 1.4} q(x,y)+q(y,x) &= 0 \\
   \label{eq 1.5} q(x,x) &=0.
  \end{align}
  
Write $q(x,y) = \sum_{i,j} c_{i,j} x^i y^j$. By the antisymmetry condition \eqref{eq 1.4}, we have $c_{i,j} + c_{j,i} = 0$ for all $i,j$. Then \eqref{eq 1.5} implies that
\[\sum_{i+j=k}c_{i,j}=0\]
for all $k$. This condition is vacuous when $k$ is odd (since the summands cancel in pairs) and, when $k = 2n$ is even, it forces $c_{n,n} = 0$. Hence, $q$ can be written as
\[q(x,y)=\sum_{i<j}c_{i,j}(x^iy^j-y^ix^j)=s(x,y)-s(y,x),\]
where we define $s(x,y)\colonequals \sum_{i<j}c_{i,j}x^iy^j$. 

Let $p_0 \in R[x,y,z]$ be such that $(p_0, q) = \gamma(s)$. Then $(p, q)$ and $(p - p_0, 0)$ define the same cohomology class. It follows that every element of $\textrm{H}^3_s(\mathbb{G}_{a,R}, \mathbb{G}_{a,R})$ can be represented by a pair of the form $(p, 0)$, for some $p \in R[x,y,z]$. Any such representative $p$ satisfies the simplified relations
\begin{align} 
   \label{eq 2.1} p(x+y,z,w)+p(x,y,z+w) &= p(y,z,w)+p(x,y+z,w)+p(x,y,z)\tag{1$^\prime$} \\
   \label{eq 2.2} p(x,y,z)+p(z,x,y) &= p(x,z,y) \tag{2$^\prime$}\\
   \label{eq 2.3} p(x,y,z)+p(y,z,x) &= p(y,x,z). \tag{3$^\prime$}
  \end{align}

Now suppose that the cohomology class of $(p,0)$ lies in the kernel of $\eta$. That is, there exists a polynomial $r \in R[x,y]$ such that
\[p(x,y,z)= r(x+y,z)+r(x,y)-r(x,y+z)-r(y,z).\]
Let $A(x,y)$ be the antissymetrization $r(x,y)-r(y,x)$ of $r$. Then equations \eqref{eq 2.2} and \eqref{eq 2.3} imply that $A$ is additive in each variable. Since $R$ is torsion-free as an abelian group, any bivariate polynomial over $R$ that is additive in both variables must be of the form $cxy$ for some $c \in R$. But such a polynomial is antisymmetric only if $c = 0$. Hence, $A = 0$ and $r$ is symmetric.

We conclude that $(p,0) = \gamma(r)$, so the class of $(p,0)$ is trivial in cohomology. This shows that the map $\eta$ is injective. In particular, $\textrm{H}^3_s(\mathbb{G}_{a,R}, \mathbb{G}_{a,R})$ vanishes whenever $\textrm{H}^3_0(\mathbb{G}_{a,R}, \mathbb{G}_{a,R})$ does. The latter group is known to vanish when $R$ is a $\mathbb{Q}$-algebra, by \cite[Thm.\ 1]{polynomialcocycles}.
\end{proof}

The same strategy used in the proof of Corollary~\ref{H2} applies here to show that Proposition~\ref{H3 Ga} yields the following multiplicative counterpart.

\begin{corollary}\label{H3}
Let $R$ be a ring. If either $R$ is reduced or a $\mathbb{Q}$-algebra, then both cohomology groups $\normalfont\textrm{H}^3_s(\mathbb{G}_{a,R}, \mathbb{G}_{m,R})$ and $\normalfont\textrm{H}^3_0(\mathbb{G}_{a,R}, \mathbb{G}_{m,R})$ vanish.
\end{corollary}

\nocite{*}
\printbibliography

\textsc{Gabriel Ribeiro, Department of Mathematics, ETH Zurich, 8092 Zurich, Switzerland}

\textit{Email address:} \href{mailto:gabriel.ribeiro@math.ethz.ch}{gabriel.ribeiro@math.ethz.ch}

\vspace{1em}

\textsc{Zev Rosengarten, Einstein Institute of Mathematics, The Hebrew University of Jerusalem, Edmond J. Safra Campus, 91904, Jerusalem, Israel}

\textit{Email address:} \href{mailto:zevrosengarten@gmail.com}{zevrosengarten@gmail.com}

\end{document}